\documentclass[anz]{cupaus}

\usepackage[T1]{fontenc}
\usepackage{amsfonts}
\usepackage{amssymb}
\usepackage{latexsym}
\usepackage{color}
\usepackage{graphicx}
\usepackage{algorithmic}
\usepackage{eurosym}

\usepackage{algorithmic}
\usepackage{amsmath}
\usepackage{amssymb}
\usepackage{bbm}
\usepackage{eurosym}
\usepackage[T1]{fontenc}
\usepackage{upgreek}
\usepackage{mathtools}

\usepackage[numbers,sort&compress]{natbib}

\newtheorem{problem}{Problem}
\newtheorem{example}{\textit{Example}}

\newtheorem{remark}{Remark}

\newcommand{\bft}{\mbox{$\mbox{\boldmath $t$}$}} 
\newcommand{\bfu}{\mbox{$\mbox{\boldmath $u$}$}} 
\newcommand{\bfw}{\mbox{$\mbox{\boldmath $w$}$}} 

\newcommand{\sbfw}{\mbox{\boldmath \scriptsize $w$}} 

\newcommand{\bfJ}{\mbox{$\mbox{\boldmath $J$}$}}

\newcommand{\bfU}{\mbox{$\mbox{\boldmath $U$}$}} 
\newcommand{\bfV}{\mbox{$\mbox{\boldmath $V$}$}}

\newcommand{\bfzero}{\mbox{$\mbox{\boldmath $0$}$}} 

\newcommand{\bdot}{\scalebox{1.5}{.}}

\begin{document}

\verso{Optimal driving strategies for a fleet of trains \texttt{cupaus.cls}?}
\recto{Howlett, Kapsis and Pudney}

\title{Optimal driving strategies for a fleet of trains subject to total energy consumption constraints on predetermined time intervals}

\cauthormark
\author[1]{Phil Howlett}
\author[2]{Maria Kapsis}
\author[3]{Peter Pudney}

\address[1]{Scheduling and Control Group, UniSA STEM, University of South Australia, South Australia, 5095, Australia \email[1]{phil.howlett@unisa.edu.au}}
\address[2]{Scheduling and Control Group, UniSA STEM, University of South Australia, South Australia, 5095, Australia \email[2]{maria.kapsis@mymail.unisa.edu.au}}
\address[3]{Scheduling and Control Group, UniSA STEM, University of South Australia, South Australia, 5095, Australia \email[2]{peter.pudney@unisa.edu.au}}
      
\pages{1}{26}
      
\begin{abstract}
In order to manage electricity transmission and distribution it is now common practice for system operators to offer financial incentives that encourage large consumers to reduce energy usage during designated peak demand periods.  For train operators on large rail networks it may be profitable\textemdash with selected individual journeys\textemdash to reduce energy usage during peak times and increase energy usage at other times rather than simply minimizing overall energy consumption.  We will use classical methods of constrained optimization to find optimal driving strategies for a fleet of trains subject to limits on total energy consumption during specified intermediate time intervals but with no change to individual journey times.  The proposed strategies can be used by a large rail organisation to reduce overall operating costs with only minimal disruption to existing schedules and with no changes to important departure and arrival times.
\end{abstract}

\keywords[MSC 2020]{\textit{49K15, 49M15, 93C95}}

\keywords[\textit{Key words and phrases}]{Train control, optimal driving strategies, energy consumption constraints}

\maketitle

\section{Introduction}
\label{s:int}

In order to improve management of supply and demand for electrical energy it is now common practice for system operators to offer financial incentives that encourage high-volume users to reduce consumption during times of peak demand.  In this paper we will show that large rail operators can use these incentive payments to reduce overall energy costs without serious disruption to existing schedules by reducing optimal driving speeds during periods of high demand and increasing optimal driving speeds at other times.  If individual trains are already using driving strategies that minimize overall energy consumption for the scheduled journey then the suggested changes will increase overall energy consumption.  Nevertheless we will show that if we find modified optimal driving strategies that minimize overall energy consumption subject to energy consumption constraints on the designated high-demand time intervals then suitable incentive payments from the electricity company to the rail operator can more than compensate for the overall increase in energy consumption.  Thus the net  journey cost decreases even though total energy consumption has increased.

\subsection{Application}
\label{ss:mot}

In recent years large rail organisations such as the French National Rail Operator \textit{Soci\'{e}t\'{e} Nationale des Chemins de fer Fran\c{c}ais (SNCF)} have substantially reduced overall energy costs using driving strategies that minimize energy consumption but preserve scheduled journey times.  Drivers on the premier high-speed passenger service \textit{Train {\`a} Grande Vitesse (TGV)} currently use an on-board Driver Advice System (DAS) in the form of a tablet app \textit{Opticonduite} based on the \textit{Energymiser}\textsuperscript{\textregistered} system\footnote{\textit{Energymiser}\textsuperscript{\textregistered} was originally developed by the Scheduling and Control Group at the University of South Australia for \textit{TTG Transportation Technology}.  The rights are now owned by \textit{Trapeze Rail}.} to find driving strategies that minimize energy consumption and assist with on-time arrival.  The strategies proposed in this paper can easily be implemented using the \textit{Energymiser}\textsuperscript{\textregistered} technology.  A preliminary study for \textit{SNCF} into the feasibility of using modified driving strategies to reduce energy consumption at peak demand times was reported in a 2022 paper by Pam \textit{et al.}~\cite{pam1}.  A recent paper by Howlett \textit{et al.}~\cite{how13} found necessary and sufficient conditions for optimal switching in these modified strategies.  A subsequent paper by Howlett \textit{et al.}~\cite{how14} shows that these modified strategies can be further improved by considering a wider class of strategies.

\subsection{Preliminary problem description}
\label{iss:ppd}

Suppose that a large rail organisation operates a fleet of $m$ identical trains $\{ {\mathfrak T}_j\}_{j=1}^m$ on different tracks and that all scheduled journeys begin at time $t = 0$ and finish at time $t = T$.  Suppose further that there are fixed times $\{ t_j \}_{j=0}^{n+1}$ with $0 = t_0 < t_1 < \cdots < t_n < t_{n+1} = T$ and that the operator wishes to limit total energy consumption by the entire fleet on the interval $(t_k, t_{k+1})$ to a specified amount for each $k =0,\ldots,n$.  Find driving strategies for each train that preserve the scheduled journey time and minimize energy consumption for each individual train subject to satisfying a total energy consumption constraint for the entire fleet on each interval $(t_k, t_{k+1})$.  

\begin{remark}
\label{rem1}
The assumptions that the trains are identical and that all journeys begin at $t=0$ and end at $t = T$ are not necessary for the validity of our analysis.  We could simply assume that the equations of motion for each train take the same form but the numerical values of the parameters are different and that train ${\mathfrak T}_j$ begins at time $t_{j,0} \leq t_0$ and finishes at time $T_j \geq T$.  If we do this then the analysis remains essentially the same but the notation becomes more complicated. 
\end{remark}   

\subsection{Summary of the main results}
\label{iss:smr}

We show that the driving strategy for each individual train ${\mathfrak T}_j$ is determined by a single unconstrained optimal driving speed $V_j$ on those prescribed intervals with no active energy consumption constraint and a collection of separate constrained optimal driving speeds $V_{j,k} < V_j$ on each interval $(t_k,t_{k+1})$ with an active energy consumption constraint.  We also show that there is a vector $\bfw = (w_1,\ldots,w_n) \geq \bfzero$ of weight factors associated with each the interval $(t_k, t_{k+1})$ and that $w_k > 0$ for those intervals with an active constraint and $w_k = 0$ otherwise.  For each train ${\mathfrak T}_j$ the constrained optimal driving speeds $V_{j,k}$ are related to the unconstrained optimal driving speed $V_j$ by the relationship
$$
1 + w_k = \varphi^{\, \prime}(V_j)/\varphi^{\, \prime}(V_{j,k})
$$
where $\varphi(V)$ is the power required to drive at constant speed $V$ over a non-zero time interval.  For the given fixed journey time and a fixed value of $\bfw = (w_1,\ldots,w_n)$ the optimal unconstrained driving speed $V_j$ must be adjusted so that the distance travelled by train ${\mathfrak T}_j$ in the allocated time is equal to the specified length of the journey. The constrained optimal driving speeds $V_{j,k-1}, V_{j,k}, V_{j,k+1}$ determine the energy consumed by train ${\mathfrak T}_j$ on the interval $(t_k,t_{k+1})$.  We will show that the individual speed profiles for each train take the same general form as the strategies described recently by Howlett \textit{et al.}~\cite{how14} for the corresponding single train problem.  We will apply our new results to archetypal examples for a fleet of three trains with one active energy consumption constraint, and a fleet of five trains with three active energy consumption constraints on consecutive intervals.    

\subsection{Structure of the paper}
\label{iss:sp}

The remainder of the paper is structured as follows.  Section~\ref{s:not} summarises the notation.  Section~\ref{s:lr} is the literature review.  The review consists of a preamble which provides a general survey of the work on optimal train control and a subsection which reviews previous work on optimal driving strategies for trains subject to energy consumption constraints on pre-determined intermediate time intervals.  This section also contains a specific example which illustrates the optimal speed profiles for a single train with one active energy constraint and a complete range of feasible restrictions.  Section~\ref{s:btm} describes the basic train model and outlines how this model has been used in practice to find optimal driving strategies for real trains.  In Section~\ref{s:ocp} we formulate and solve the mathematical control problem for a fleet of trains where total energy consumption is restricted on multiple predetermined intermediate time intervals.  The Pontryagin principle is applied to show that the only allowable controls in an optimal strategy for train ${\mathfrak T}_j$ are phases of \textit{maximum acceleration}, \textit{speedhold} at the unconstrained optimal driving speed $V_j = V_{j,\ell}$ on intervals $(t_{\ell}, t_{\ell+1})$ with no active energy constraint, \textit{speedhold} at speed $V_{j,k}$ on intervals $(t_k,t_{k+1})$ with an active energy constraint, \textit{coast} and \textit{maximum brake}.  We also study the evolution of a modified adjoint variable $\eta$ and use the continuity of $\eta$ to find optimal switching points between the various phases of optimal control for each train.  We also show how the theory can be applied in specific cases to find optimal speed profiles for a fleet of three trains with a single active energy constraint and for a fleet of five trains with three consecutive active intermediate energy constraints.  In Section~\ref{s:con} we draw some brief conclusions.

\section{Notation}
\label{s:not}

All quantities are measured in SI units: time in seconds (s), distance in metres (m), mass in kilograms (kg) and energy in joules (J). The following notation and terminology is used for train ${\mathfrak T}_j$ throughout.  We write
\begin{itemize}
\item independent variable: time $t$,
\item dependent variables: position $x_j$, speed $v_j$, 
\item journey distance $X_j$, journey time $T$,
\item starting time $0$, finishing time $T$,
\item applied tractive acceleration $u_{j,a}$,
\item upper bound on tractive acceleration: $H_a(v)$,
\item applied braking force per unit mass: $u_{j,b}$,
\item upper bound on braking force per unit mass: $H_b(v)$,
\item frictional force per unit mass: $r(v)$,
\item auxiliary functions: $\varphi(v) = vr(v)$, $\psi(v) = v^2r^{\, \prime}(v)$,
\item ${\mathcal A} \subset \{1,\ldots,n\}$ is the index set for $(t_k,t_{k+1})$ with an active energy constraint,
\item unconstrained optimal driving speed: $V_j = V_{j,\ell}$ when $\ell \notin {\mathcal A}$, and
\item constrained optimal driving speed: $V_{j,k}$ when $k \in {\mathcal A}$.
\end{itemize}
The \textit{intermediate energy consumption constraints} will normally be referred to as simply \textit{the energy constraints}. 

\section{Literature review}
\label{s:lr}

We will not attempt a comprehensive review of the well-established and clearly-defined theoretical work on optimal control of a single train.  For details of the underlying theory we refer readers to the key papers on continuous control \cite{alb4, alb9, alb10, bar2, how7, how12, khm1, liu1}.  The frequently-cited references on discrete control \cite{che1, how3, how6, how8} are an important component of the literature but are less relevant to the work in this paper.  The early papers on continuous control \cite{asn1, how1} contain some important insights but these results have generally been superseded.  Other significant studies in the Russian literature \cite{bar1, gol2} may be difficult to obtain and the relevant results can usually be extracted from more recent papers.  This whole body of work relies on classical methods of constrained optimization\textemdash the Pontryagin maximum principle \cite{dub1, gir1, har1}, the Euler-Lagrange equations \cite[pp 179\textendash 190]{lue1} and the Karush-Kuhn-Tucker conditions \cite[pp 247\textendash 254]{lue1}.  These classical strategies have been implemented in real-time Driver Advice Systems (DAS's) on very fast trains where updated strategies can be recalculated on-board in a matter of a few seconds \cite{alb10}.  These systems reduce fuel costs by $5$\textendash$15\%$ and improve on-time arrival. 

There is a separate theory relating to the study of optimal driving strategies for heavy-haul freight trains where intra-train forces generated by gradient changes over the entire length of a very long train may become a significant problem.  These forces and ill-judged changes to the driving or braking forces can cause serious problems including large impulsive inertial forces which may lead to catastrophic failure of the linkages between successive wagons.  While this work is partly based on the classical theory of optimal driving the main concern is with engineering considerations where the train is modelled as a sequence of individual wagons within which collective inertial forces are transmitted via the links between successive units \cite {cho1, zhu3}. 

There is a vast body of broadly-based work on the development of efficient timetables.  The objectives include improved efficiencies obtained by adjusting allowed journey times \cite{gup1, lilo1, lilo2}, cost savings by choosing the best ordering of scheduled services, selection of best meeting points and energy recovery from regenerative braking \cite{bur2, cap1, dor1, lius2}, schedule recovery from disruption \cite{wanp2}, improved service to customers \cite{yan1} and train trajectory planning subject to safe separation constraints for trains travelling on the same track in the same direction \cite{wany1}.  There is also a collection of recent papers which seeks to develop a coherent theory of optimal driving subject to safe separation constraints for successive trains on the same track \cite{alb13, how11}.

A comprehensive survey covering all aspects of optimal train control can be found in recent review articles by Scheepmaker et al.~\cite{sch1} and Yin et al.~\cite{yin1}.

\subsection{Optimal driving strategies subject to energy consumption constraints}
\label{ss:odsecc}

This paper is most closely related to recent papers by Pam \textit{et al.}~\cite{pam1}, Kapsis \textit{et al.}~\cite{kap1} and Howlett \textit{et al.}~\cite{how13, how14}.  The paper by Pam \textit{et al.}~describes the results of a recent feasibility study for the French National Rail Operator \textit{Soci\'{e}t\'{e} Nationale des Chemins de fer Fran\c{c}ais} (\textit{SNCF}) into potential cost savings that could be made with minimal disruption to existing train schedules for the premier high-speed passenger service \textit{Train {\`a} Grande Vitesse (\textit{TGV})} by reducing energy consumption for selected individual trains during periods of peak electricity demand.  Kapsis \textit{et al}~\cite{kap1}  used simplified train models to provide initial insights into how a fleet of trains can be controlled to manage peak demand energy.  The constant-speed strategies are not realistic but are nevertheless a useful way to obtain indicative results in large network problems.  Howlett \textit{et al}.~\cite{how13} used classical methods of constrained optimization to find optimal switching points for the heuristic strategies that were employed in the study by Pam \textit{et al.}~\cite{pam1}.  Pam \textit{et al.}~assumed that trains would simply use a phase of \textit{coast} to move from a higher optimal driving speed to a lower optimal driving speed and a phase of \textit{maximum acceleration} to move from a lower optimal driving speed to a higher driving speed.  Howlett \textit{et al.}~showed that a strategy in this general class is optimal if and only if the switching points are located at the endpoints of the designated intervals $(t_k, t_{k+1})$.  A subsequent paper by Howlett \textit{et al.}~\cite{how14} showed that the strategies in~\cite{how13} were suboptimal and that the true optimal restricted strategies can be found by optimizing within a wider class of strategies of optimal type.  In this paper we will show that the optimal strategies proposed in~\cite{how14} for individual trains are also optimal for each train in a fleet of trains where the total energy consumption for the whole fleet is constrained on certain predetermined intermediate time intervals.

When there is a single intermediate interval $(t_1, t_2) \subset (0,T)$ on which energy consumption is constrained it was shown in~\cite{how14} that on level track the optimal strategy for a single train ${\mathfrak T}_1$ is \textit{maximum acceleration}, \textit{speedhold} at the optimal unconstrained driving speed $V_{1,0} = V_1$, \textit{maximum acceleration} to speed $v_1(t_1) = W_{1,1} > V_1$ at $t = t_1$, \textit{coast}, \textit{speedhold} at the constrained optimal driving speed $V_{1,1} < V_1$, \textit{coast} to speed $v_1(t_2) = W_{1,2} < V_{1,1}$ at $t = t_2$, \textit{maximum acceleration}, \textit{speedhold} at the optimal unconstrained driving speed $V_{1,2} = V_1$, \textit{coast} to the optimal braking speed $U_1 = V_1 - \varphi(V_1)/\varphi^{\, \prime}(V_1)$ and \textit{maximum brake}.

\begin{figure}[ht]
\begin{center}
\includegraphics[width=0.95\textwidth]{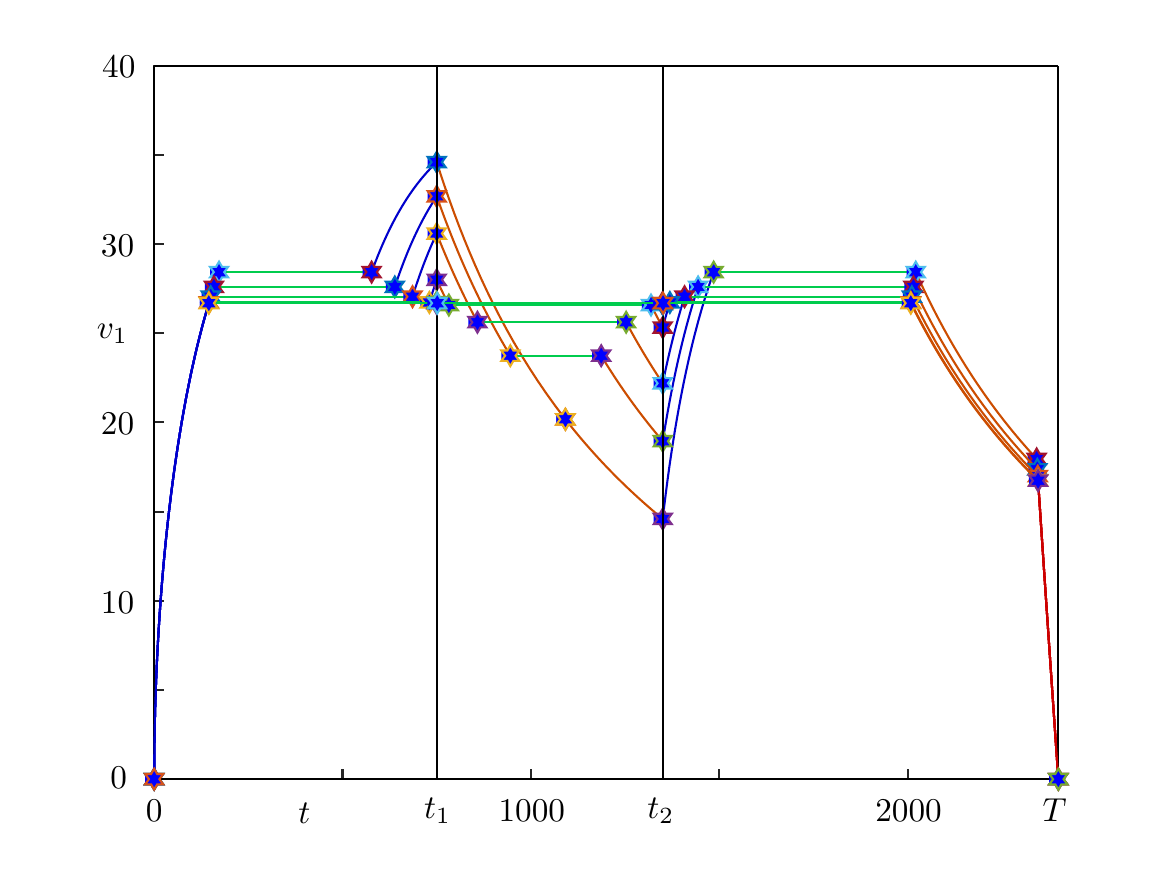}
\end{center}
\vspace{-0.5cm}
\caption{Optimal speed profiles $v_1 = v_1(t)$ in Example~\ref{ex:1}.}
\label{fig1}
\end{figure}

\begin{example}
\label{ex:1}
We consider a single train ${\mathfrak T}_1$ and use the model described in Section~\ref{s:btm} with $r(v) = r_0 + r_2v^2$ where $r_0 = 6.75 \times 10^{-3}$ and $r_2 = 5 \times 10^{-5}$,  $H_a(v) = A/v$ where $A = 3$ and $H_b(v) = 3 \times 10^{-1}$.  We assume the train has mass $M$~kg.  The journey distance is $X = 60000$ m and the journey time is $T= 2400$ s.  We used a multi-dimensional Newton iteration to find optimal speed profiles $v_1 = v_1(t)$ subject to a full range of energy consumption constraints $E_{1,1} \leq {\mathcal Q}_1 \in \{0,200 M,400 M,600 M,675 M\}$~J on the $10$\textendash minute interval $t \in (t_1,t_2) = (750, 1350)$ s.  The cost of the journey is the total energy consumption $J_1 = M \int_0^T u_{1,a}(t)v_1(t) dt$~J.  The values of the key variables are shown in TABLE~1 and the corresponding optimal speed profiles are shown in Figure~\ref{fig1}.  When ${\mathcal Q}_1 = 0$~J the train \textit{coasts} during the entire interval $t \in (t_1,t_2)$ and the total cost is $J_1 \approx (2702 M)$ J.  When ${\mathcal Q}_1 > (675 M)$~J  the intermediate energy consumption constraint becomes inactive, the total cost is $J_1 \approx (2541 M)$~J and the speed profile reverts to the classic optimal speed profile, \textit{maximum acceleration}, \textit{speedhold}, \textit{coast} and \textit{maximum brake}.

If the normal cost of energy is $c$ dollars per joule and the train mass is $M$~kg then the rail company would need to receive a reduction in their energy bill that exceeds the normal cost of the additional energy used for the modified strategy.  For this particular train the rail operator would agree to limit energy consumption on the high-demand interval to $(400 M)$~J if the electricity company agreed to pay the operator significantly more than the difference between the cost of the restricted journey \$$(2551 M)c$ and the cost of the optimal journey \$$(2541M) c$.  That is \$$(10 M)c$.  In practice this is a matter of scale.  The energy company would only do this if the reduction in energy consumption on the high-demand interval was sufficiently large.  Thus the idea would not be implemented at the level of a single train but would be feasible for a sufficiently large fleet of trains. $\hfill \Box$
\end{example}

\begin{table}[ht]
\vspace{0.25cm}
\begin{center}
\caption{Values of the key variables in Example~\ref{ex:1}.}
\begin{tabular}{|c|c|c|c|c|c|c|} \hline
\rule{0cm}{0.25cm} ${\mathcal Q}_1$ & $J_1$ & $V_1$ & $W_{1,1}$ & $V_{1,1}$ & $W_{1,2}$ & $U_1$ \\ \hline
\rule{0cm}{0.25cm} 0 & $2702 M$ & 28.43  &  34.59  &  20.17  &  14.59  &  17.95 \\ \hline
\rule{0cm}{0.25cm} $200 M$ & $2592 M$ & 27.59  & 32.67  &  23.74  & 18.85  &  17.37 \\ \hline
\rule{0cm}{0.25cm} $400 M$ & $2551 M$ & 27.04  & 30.59  & 25.62  & 22.19  & 16.98 \\ \hline
\rule{0cm}{0.25cm} $600 M$ & $2541 M$ & 26.72 & 28.00  & 26.58  & 25.32  & 16.76 \\ \hline
\rule{0cm}{0.25cm} $675 M$ & $2541 M$ & 26.68 & 26.74  & 26.68 & 26.63 & 16.73 \\ \hline
\end{tabular}
\end{center}
\end{table}

\section{The basic train model}
\label{s:btm}

To begin we note that Howlett and Pudney~\cite[Section 2.3]{how4} showed that the motion of a train with distributed mass can be reduced to the motion of a point-mass train.  Thus we restrict our attention to point-mass trains.  For further discussion of the modelling issues we refer to \cite{alb9, alb10, how4, how7, khm1, liu1}.  The motion of a train can be described using three state variables\textemdash the position $x \in [0, X]$, the elapsed journey time $t \in [0, T]$ and the speed $v \in [0, \infty)$.  The standard model uses $x$ as the independent variable with speed $v = v(x)$ and elapsed journey time $t = t(x)$ as dependent variables.  In this paper we use a non-standard model with $t$ as the independent variable and with $x = x(t)$ and $v = v(t)$ as the dependent variables.  The two models are equivalent.  We measure distance in metres (m) and time in seconds (s).  Consequently the speed $v = v(t)$ is measured in metres per second (m s$^{-1}$).  The equations of motion are
\begin{eqnarray}
\overset{\bdot}{x} & = & v \label{e:xt} \\
\overset{\bdot}{v} & = & u_a - u_b - r(v) + g(x) \label{e:vt}
\end{eqnarray}
where $(x,v) = (x(t), v(t))$ is the state vector for $t \in [0,T]$, $u_a = u_a(t) \geq 0$ is the controlled driving force per unit mass or forward acceleration and $u_b = u_b(t) \geq 0$ is the magnitude of the controlled braking force per unit mass both measured in metres per second squared (m s$^{-2}$).  We have written $\overset{\bdot}{x} = dx/dt$ and $\overset{\bdot}{v} = dv/dt$.  For a journey from $t = 0$ to $t = T$ we assume that $x(0) = 0$ and $x(T) = X$ and $v(0) = v(T) = 0$ with $v(t) > 0$ for all $t \in (0,T)$.  We assume that $u_a(t)$ and $u_b(t)$ are bounded with $0 \leq u_a(t) \leq H_a[v(t)]$ and $0 \leq u_b(t) \leq H_b[v(t)]$ for each $t \in [0, T]$ where $H_a(v), H_b(v) \in (0,\infty)$ for $v \in (0, \infty)$ are bounded monotone functions with $H_a(v) \downarrow 0$ as $v \uparrow \infty$.  The functions $H_a$ and $H_b$ define bounds for the maximum braking and driving forces per unit mass in a form that includes\textemdash as special cases\textemdash the specified bounds for a wide range of modern electric and diesel-electric locomotives.  The function $r(v)$ is a general backward force per unit mass due to frictional resistance measured in metres per second squared (m s$^{-2}$) with no specific formula assumed.  We define auxiliary functions $\varphi(v) = vr(v)$ and $\psi(v) = v^2 r^{\, \prime}(v)$ and assume only that $\varphi(v)$ is strictly convex with $\varphi(v) \geq 0$ for $v \geq 0$ and $\varphi(v)/v \uparrow \infty$ as $v \uparrow \infty$.  Note that $d[\varphi(v)/v]/dv = \psi(v)/v^2 > 0$.  See~\cite[Appendix A.3, p 409]{how6} for an expanded discussion.  It follows that both $r(v)$ and $\psi(v)$ are non-negative and strictly increasing for $v \geq 0$.  These properties capture the functional characteristics of the traditional quadratic resistance formula\textemdash the so-called Davis formula \cite{dav1}.  The function $g(x)$ is nominally the component of gravitational acceleration due to track gradient but in practice may also include additional position-dependent resistive forces.  We assume that $g(x)$ is piecewise continuous.  We also assume that energy recovered from regenerative braking is not used to drive the train.  This assumption generates a more relaxed optimal driving strategy that encourages coasting and discourages braking.  The cost of the journey is the mechanical energy consumption in joules ($\equiv$ kg m$^2$ s$^{-2}$) given by
\begin{equation}
\label{e:cost}
J = M \mbox{$\int_0^T$} u_a(t)v(t) \, dt.
\end{equation}
Suppose a particular speed $v = v_0(t)$ for $t \in [0,T]$ is prescribed.  We can use (\ref{e:xt}) to find the corresponding position $x(t)$ for all $t \in [0,T]$ and then we can use (\ref{e:vt}) to find $u_a(t)$ and $u_b(t)$.  The prescribed speed is feasible if $-H_b[v_0(t)] \leq u_a(t) -  u_b(t) \leq H_a[v_0(t)]$ for all $t \in [0,T]$.  An interval $(b,c) \subset [0, X]$ is said to be steep uphill at speed $v = W$ if
\begin{equation}
\label{e:suw}
H_a(W) - r(W) + g(x) < 0
\end{equation}
for all $x \in (b,c)$.  An interval $(b, c)$ is said to be steep downhill at speed $v = W$ if
\begin{equation}
\label{e:sdw}
- r(W) + g(x) > 0
\end{equation}
for all $x \in (b,c)$.  Whether or not a track is classified as steep depends on the desired speed $W$ of the train.

\subsection{Optimal train control in practice}
\label{ss:otcp}

The use of optimal speed profiles reduces fuel costs but as a journey evolves and circumstances change it is necessary to continually recalculate the optimal profile.  Specialized numerical algorithms developed by the Scheduling and Control Group (SCG) at the University of South Australia for \textit{TTG Transportation Technology} and now owned by \textit{Trapeze Rail} are key components of the Driver Advice System (DAS) known in the industry as \textit{Energymiser}\textsuperscript{\textregistered}.  The system uses numerical calculation routines based on the classical theory of optimal control for a single train and routinely computes optimal speed profiles for freight and passenger trains on non-level tracks with steep gradients for journey segments of more than 100 kms in a few seconds.  \textit{Energymiser}\textsuperscript{\textregistered} is used by major rail operators in the United Kingdom to provide on-board advice in real time to train drivers about energy-efficient driving strategies and the corresponding speed profiles.  Since 2017 \textit{Energymiser}\textsuperscript{\textregistered} has been deployed by the French National Railway Company {\em Soci\'{e}t\'{e} Nationale des Chemins de fer Fran\c{c}ais} (\textit{SNCF}) as a tablet app by all drivers on all services including the famed high-speed electric rail service {\em Train \`{a} Grande Vitesse} (\textit{TGV}).  The system, known as {\em Opti-conduite} in France, was the 2017 winner of {\em Enterprise of the Year} awarded by Directors of Energy Departments in France.  See \cite[Section 9]{alb10} for details of the \textit{TGV} application.  The system reduces fuel consumption by between 5\% and 15\% and improves on-time arrival.  The cost of energy (electricity) for \textit{SNCF} exceeded \euro 1,300,000,000 in 2021.  More information about {\em Energymiser}\textsuperscript{\textregistered} can be obtained from the \textit{Trapeze Rail} website\footnote{See https://trapezegroup.com.au/resources/on-time-passenger-rail-operators-trapeze-das/}.  General methods of computational control are not well suited to real-time calculation of optimal driving strategies and cannot compete in practice with purpose-built systems \cite{sch1}.

\section{The optimal control problem}
\label{s:ocp}

Let $\{ t_k \}_{k=0}^{n+1}$ be a strictly increasing sequence of fixed times with $t_0 = 0$ and $t_{n+1} = T$ and let $\{{\mathfrak T}_j\}_{j=1}^m$ be a fleet of $m$ trains.  Suppose that each train ${\mathfrak T}_j$ starts at time $t = 0$ and finishes at time $t = T$.  We assume that the journeys are independent parallel journeys on level tracks and that train ${\mathfrak T}_j$ travels from $x_j = 0$ to $x_j = X_j$.  We suppose that energy consumption across the entire fleet is limited to a specified total amount ${\mathcal Q}_k$ on each interval $t \in (t_k,t_{k+1})$ for $k=0,1,\ldots,n$.  It is convenient to assume that the trains are identical.  In this case we may assume without loss of generality that $M = 1$.  We assume that all specified journeys are feasible.  The results remain valid for non-identical trains with different journey times provided each train ${\mathfrak T}_j$ starts at some time $t_{j,0} \leq t_0$ and finishes at some time $T_j \geq T$.  For non-identical trains the masses of the individual trains must be used in the formulation of the cost function.

\begin{problem}
\label{p:1} 
For each $j=1,\ldots,m$ find the controls $u_{j,a} = u_{j,a}(t)$ and $u_{j,b} = u_{j,b}(t)$ and the associated functions $x_j = x_j(t)$ and $v_j = v_j(t)$ satisfying $(\ref{e:xt})$ and $(\ref{e:vt})$ that minimize total energy consumption
\begin{equation}
\label{e:cp1}
J(\bfu) = \mbox{$\sum_{j=1}^m \int_0^T$} u_{j,a}(t) v_j(t) dt
\end{equation}
subject to the control constraints $u_{j,a} \in [0, H_a(v)]$ and $u_{j,b} \in [0, H_b(v)]$, the energy consumption constraints
\begin{equation}
\label{e:eccp1}
E_k(\bfu) - {\mathcal Q}_k \leq 0
\end{equation}
where $E_k(\bfu) = \mbox{$\sum_{j=1}^m \int_{t_k}^{t_{k+1}}$} u_{j,a}(t)v_j(t) dt$ is the total energy consumption when $t \in (t_k, t_{k+1})$ for each $k = 0,1,\ldots,n$, and the boundary conditions $(x_j(0), v_j(0)) = (0,0)$ and $(x_j(T),v_j(T)) = (X_j,0)$.
\end{problem}

We apply the Pontryagin principle and follow the arguments used by Howlett \textit{et al.}~\cite{how14} but generalize the analysis to a fleet of trains rather than a single train.  Readers are referred to \cite{how14} for a more detailed presentation and for additional commentary on the methodology.

\subsection{The weighted cost function}
\label{ss:wcf}

The weighted cost function is
\begin{eqnarray}
\label{e:wcf}
{\mathcal J} & = & \mbox{$\sum_{j=1}^m$} \mbox{$\int_0^T$} u_{a,j}(t) v_j(t) dt + \mbox{$\sum_{k=0}^n$} w_k \left[ E_k(\bfu) - {\mathcal Q}_k \right] \nonumber \\
& = & \mbox{$\sum_{j=1}^m$} \mbox{$\int_0^T$} \left[ \rule{0cm}{0.4cm} 1 + \mbox{$\sum_{k=0}^n$} w_k {\mathbbm 1}_k(t) \right] u_{a,j}(t) v_j(t) dt \nonumber \\
& & \hspace{3cm} - \mbox{$\sum_{k=0}^n$} w_k {\mathcal Q}_k
\end{eqnarray}
where
$$
{\mathbbm 1}_k(t) = \left\{ \begin{array}{ll}
1 &\mbox{for}\ t \in (t_k,t_{k+1}) \\
1/2 &\mbox{for}\ t = t_k, t_{k+1} \\
0 &\mbox{otherwise}. \end{array} \right.
$$
is the indicator function for the interval $[t_k, t_{k+1}]$ and $\bfw = (w_k)_{k=0}^n$ is a vector of non-negative weights.  If $w_k > 0$ then the minimization ensures that the constraint (\ref{e:eccp1}) is active. 

\subsection{Hamiltonian function}
\label{ss:hf}

The Hamiltonian function is
\begin{eqnarray}
\label{e:hf}
{\mathcal H}_{\sbfw} & = & \mbox{$\sum_{j=1}^m$} \left\{ \rule{0cm}{0.4cm} (-1) \left[ 1 + \mbox{$\sum_{k=0}^n$} w_k {\mathbbm 1}_k \right] u_{a,j} v_j \right. \nonumber \\
& &  \left. \rule{0cm}{0.4cm} + \lambda_j v_j + \mu_j [u_{a,j} - u_{b,j} - r(v_j) + g(x_j)] \right\}
\end{eqnarray}

\subsection{Lagrangian function}
\label{ss:lf}

The Lagrangian function is
\begin{eqnarray}
\label{e:lf}
{\mathcal L}_{\sbfw} & = & {\mathcal H}_{\sbfw} + \mbox{$\sum_{j=1}^m$} \left\{ \rule{0cm}{0.4cm} \alpha_j u_{a,j} + \beta_j[H_a(v_j) - u_{a,j}] \right. \nonumber \\
& & \hspace{1cm} \left. \rule{0cm}{0.4cm} + \gamma_j u_{b,j} + \delta_j[ H_b(v_j) - u_{b,j}] \right\}
\end{eqnarray}
where $\alpha_j$, $\beta_j$, $\gamma_j$, $\delta_j$ are non-negative Lagrange multipliers.

\subsection{Adjoint equations}
\label{ss:ae}

The adjoint equations are given by
\begin{eqnarray}
\overset{\bdot}{\lambda_j} & = & - \frac{\partial {\mathcal L}_{\sbfw}}{\partial x_j} \nonumber \\
& = & - \mu_j g^{\, \prime}(x_j). \label{e:ae1} \\
& & \nonumber  \\
\overset{\bdot}{\mu_j} & = & - \frac{\partial {\mathcal L}_{\sbfw}}{\partial v_j} \nonumber \\
& = & (-1) \left[ 1 + \mbox{$\sum_{k=0}^n$} w_k {\mathbbm 1}_k \right] u_{a,j} - \lambda_j + \mu_j r^{\, \prime}(v_j) \nonumber \\
& & \hspace{1cm} - \beta_j H_a^{\, \prime}(v_j) - \delta_j H_a^{\, \prime}(v_j). \label{e:ae2}
\end{eqnarray}
for each $j =1,\ldots,m$.

\subsection{Karush\textendash Kuhn\textendash Tucker equations}
\label{ss:kkt}

The Karush\textendash Kuhn\textendash Tucker (KKT) equations are
\begin{eqnarray*}
\frac{\partial {\mathcal L}_{\sbfw}}{\partial u_{a,j}} = 0 & \Leftrightarrow & (-1)(1 + w_k{\mathbbm 1}_k)v_j + \mu_j + \alpha_j - \beta_j = 0, \\
\frac{\partial {\mathcal L}_{\sbfw}}{\partial u_{b,j}} = 0 & \Leftrightarrow & - \mu_j + \gamma_j - \delta_j = 0.
\end{eqnarray*}

\subsection{Complementary slackness conditions}
\label{ocp:cs}

The complementary slackness conditions are
\begin{eqnarray*}
\alpha_j u_{a, j} & = & 0,  \\
\beta_j[ H_a(v_j) - u_{a,j}] & = & 0, \\
\gamma_ju_{b, j} & = & 0, \\
\delta_j[H_b(v_j) - u_{b, j}] & = & 0.
\end{eqnarray*}

\subsection{The optimal controls}
\label{ss:oc}

For train ${\mathfrak T}_j$ the optimal control is determined by the relative values of the state variable $v_j$ and the adjoint variable $\mu_j$.   We assume $t \in (t_k, t_{k+1})$ for some fixed $k \in \{0,1,\ldots,n\}$.

\vspace{0.2cm}
\textbf{Case 1: $\mu_j > (1 + w_k)v_j$.}\, $\partial {\mathcal L}_{\sbfw}/\partial u_{a,j} = 0 \Rightarrow \mu_j - (1 + w_k)v_j = \beta_j - \alpha_j \Rightarrow (\alpha_j, \beta_j) = ( 0, \mu_j - (1+ w_k)v_j)$.  Therefore $\beta_j(H_a(v_j) - u_{a,j}) = 0 \Rightarrow u_{a,j} = H_a(v_j)$.  $\partial {\mathcal L}_{\sbfw}/\partial u_{b,j} = 0 \Rightarrow \mu_j = \gamma_j - \delta_j \Rightarrow (\gamma_j, \delta_j) = (\mu_j, 0)$.  Hence $\gamma_j u_{b,j} = 0 \Rightarrow u_{b,j} = 0$.  The optimal control is {\em maximum acceleration}. 

\vspace{0.2cm}
\textbf{Case 2: $\mu_j = (1 + w_k)v_j$.}\,  $\partial {\mathcal L}_{\sbfw}/\partial u_{a,j} = 0 \Rightarrow 0 = \beta_j - \alpha_j \Rightarrow (\alpha_j, \beta_j) = ( 0, 0)$.  $\partial {\mathcal L}_{\sbfw}/\partial u_{b,j} = 0 \Rightarrow \mu_j = \gamma_j - \delta_j \Rightarrow (\gamma_j, \delta_j) = (\mu_j, 0)$  Therefore $\gamma_j u_{b,j} = 0 \Rightarrow u_{b,j} = 0$.  The first adjoint equation (\ref{e:ae1}) becomes
$$
\overset{\bdot}{\lambda_j} = (-1)(1 + w_k)v_j g^{\, \prime}(x_j) = (-1)(1 + w_k)g^{\, \prime}(x_j) \overset{\bdot}{x_j}
$$
which gives $\lambda_j = - (1 + w_k)g(x_j) + \lambda_{j,k}$ for some constant $\lambda_{j,k} \in {\mathbb R}$.  Now $\mu_j = (1 + w_k)v_j$ for $t \in (t_k, t_{k+1})$ implies
$$
\overset{\bdot}{\mu_j} = (1 + w_k)\overset{\bdot}{v_j}
$$
from which we can use (\ref{e:vt}) and (\ref{e:ae2}) to deduce that $(-1)(1 + w_k)u_{a,j} - [ - (1+ w_k)g(x_j) + \lambda_{j,k}] + [(1 + w_k)v_j] r^{\, \prime}(v_j) = (1 + w_k)[u_{a,j} - r(v_j) + g(x_j)] \Rightarrow (1 + w_k)\varphi^{\, \prime}(v_j) = \lambda_{j,k} \Rightarrow v_j = V_{j,k}$ where the constant $V_{j,k} > 0$ is the unique solution to $(1 + w_k)\varphi^{\, \prime}(v_j) = \lambda_{j,k}$.  The optimal control is {\em speedhold} at speed $v_j = V_{j,k}$ using partial acceleration $u_{a,j} = r(V_{j,k}) - g(x_j)$ provided $0 \leq r(V_{j,k}) - g(x_j) \leq H_a(v_j)$.

\vspace{0.2cm}
\textbf{Case 3: $0 < \mu_j < (1 + w_k)v_j$.}\,  $\partial {\mathcal L}_{\sbfw}/\partial u_{a,j} = 0 \Rightarrow (1 + w_k)v_j - \mu_j = \alpha_j - \beta_j \Rightarrow (\alpha_j, \beta_j) = ((1 + w_k)v_j - \mu_j, 0)$.  Hence $\alpha_j u_{a,j} = 0 \Rightarrow u_{a,j} = 0$.  $\partial {\mathcal L}_{\sbfw}/\partial u_{b,j} = 0 \Rightarrow \mu_j = \gamma_j - \delta_j \Rightarrow (\gamma_j, \delta_j) = (\mu_j, 0)$.  Therefore $\gamma_j u_{b,j} = 0 \Rightarrow u_{b,j} = 0$.  The optimal control is {\em coast}.

\vspace{0.2cm}
\textbf{Case 4: $\mu = 0$.}\,  $\partial {\mathcal L}_{\sbfw}/\partial u_{a,j} = 0 \Rightarrow (1 + w_k)v_j = \alpha_j - \beta_j \Rightarrow (\alpha_j, \beta_j) = ((1 + w_k)v_j, 0)$.  Thus $\alpha_j u_{a,j} = 0 \Rightarrow u_{a,j} = 0$.  $\partial {\mathcal L}_{\sbfw}/\partial u_{b,j} = 0 \Rightarrow 0 = \gamma_j - \delta_j \Rightarrow (\gamma_j, \delta_j) = (0, 0)$.  Therefore $\overset{\bdot}{\mu_j} = 0 \Rightarrow \lambda_{j,k} = 0 \Rightarrow \varphi^{\, \prime}(V_{j,k}) = 0$.  This cannot occur in a realistic problem.

\vspace{0.2cm}
\textbf{Case 5: $\mu < 0$.}\  $\partial {\mathcal L}_{\sbfw}/\partial u_{a,j} = 0 \Rightarrow (1 + w_k)v_j + \vert \mu_j \vert = \alpha_j - \beta_j \Rightarrow (\alpha_j, \beta_j) = ((1 + w_k)v_j + \vert \mu_j \vert, 0)$.  Hence $\alpha_j u_{a,j} = 0 \Rightarrow u_{a,j} = 0$.  $\partial {\mathcal L}_{\sbfw}/\partial u_{b,j} = 0 \Rightarrow \vert \mu_j \vert  = \delta_j - \gamma_j \Rightarrow (\gamma_j, \delta_j) = (0, \vert \mu_j \vert)$.  Thus $\delta_j(H_b(v_j) - u_{b,j}) = 0 \Rightarrow u_{b,j} = H_b(v_j)$.  The optimal control is {\em maximum brake}.

It is clear that the optimal driving speed $V_{j,k}$ for train ${\mathfrak T}_j$ on the interval $t \in (t_k, t_{k+1})$ depends on the value $w_k \geq 0$.  If $w_k > 0$ a cost penalty can be avoided if we set 
$$
\mbox{$\sum_{j=1}^m$} \mbox{$\int_{t_{k-1}}^{t_k}$} u_{a,j}(t)v_j(t) dt - {\mathcal Q}_k = 0.
$$
In this case we say that the constraint is active.  Define ${\mathcal A} = \{ k \in \{0,\ldots,n\} \mid w_k > 0\}$ as the set of active constraints.  For each $\ell \in \{0,\ldots,n\} \setminus {\mathcal A}$ we have $w_{\ell} = 0$ and we say that the corresponding constraint is inactive.

\subsection{Evolution of the adjoint variable on level track}
\label{ss:eavlt}

Assume the track is level with $g(x) \equiv 0$.  In this case the first adjoint equation (\ref{e:ae1}) reduces to
$$
\overset{\bdot}{\lambda_j} = 0 \Longrightarrow \lambda_j(t) = \kappa_j
$$
for some constant $\kappa_j$ and all $t \in [0,T]$.  When $\ell \notin {\mathcal A}$ we have $w_{\ell} = 0$ and $V_{j,\ell}$ is the unique solution to the equation
$$
\varphi^{\, \prime}(v_{j,\ell}) = \kappa_j
$$
for each $j=1,\ldots,m$.  It follows that we must have $\kappa_j > 0$ and that $V_{j,\ell} = V_j$ is independent of $\ell$.  We call $V_j$ the unconstrained optimal driving speed.  When $k \in {\mathcal A}$ we have $w_k > 0$ and $V_{j,k}$ is the unique solution to the equation
$$
(1 + w_k) \varphi^{\, \prime}(v_{j,k}) = \kappa_j
$$
for each $j =1,\ldots,m$.  We call $V_{j,k}$ the constrained optimal driving speed on $(t_k,t_{k+1})$.  It follows that
$$
1 + w_k = \varphi^{\, \prime}(V_j)/\varphi^{\, \prime}(V_{j,k})
$$
and hence that $V_{j,k} < V_j$.  We can also solve (\ref{e:ae2}) to find $\mu_j(t)$ but it is more convenient to define a modified adjoint variable $\eta_j = \mu_j/v_j$ for each $j \in \{1,\ldots,m\}$. We will also define the tangent function to the graph $y = \varphi(v)$ at point $v = V$ by the formula
$$
L_V(v) = \varphi(V) + \varphi^{\, \prime}(V)(v - V)
$$
for all $v \geq 0$.  The strict convexity of $\varphi(v)$ means that $\varphi(v) > L_V(v)$ for all $v \neq V$.  The optimal driving strategy for ${\mathfrak T}_j$ is uniquely defined by the unconstrained optimal driving speed $V_j = V_{\ell,j}$ on the intervals $(t_{\ell}, t_{\ell + 1})$ where $\ell \notin {\mathcal A}$ and by the constrained optimal driving speeds $V_{j,k}$ on the intervals $(t_k, t_{k+1})$ where $k \in {\mathcal A}$.  On the intervals $(t_{\ell}, t_{\ell+1})$ where $\ell \notin {\mathcal A}$ the modified adjoint variable is given by
\begin{equation}
\label{e:etaellma}
\eta_{j,a} =  \frac{v_jH_a(v_j) - L_{V_j}(v_j)}{v_jH_a(v_j) - \varphi(v_j)}
\end{equation}
for a phase of \textit{maximum acceleration}, $\eta_{j,h} = 1$ for the \textit{speedhold} phase at speed $V_j$ and by
\begin{equation}
\label{e:etaellc}
\eta_{j,c} =  \frac{L_{V_j}(v_j)}{\varphi(v_j)}
\end{equation}
for a phase of \textit{coast}.  On the intervals $(t_k,t_{k+1})$ where $k \in {\mathcal A}$ the modified adjoint variable is given by
\begin{equation}
\label{e:etakma}
\eta_{j,a} =  (1 + w_k) \frac{v_jH_a(v_j) - L_{V_{j,k}}(v_j)}{v_jH_a(v_j) - \varphi(v_j)}
\end{equation}
for a phase of \textit{maximum acceleration}, $\eta_{j,h} = 1 + w_k$ for the \textit{speedhold} phase at speed $V_{j,k}$ and by
\begin{equation}
\label{e:etakc}
\eta_{j,c} =  (1 + w_k) \frac{L_{V_{j,k}}(v_j)}{\varphi(v_j)}
\end{equation}
for a phase of \textit{coast}.  For detailed derivations of (\ref{e:etaellma}), (\ref{e:etaellc}), (\ref{e:etakma}) and (\ref{e:etakc}) we refer readers to \cite{how14}.

\subsection{The optimal strategy for \boldmath${\mathfrak T}_j$ on level track}
\label{ss:oslt}

We can now give a categorical description of the optimal restricted strategy for train ${\mathfrak T}_j$ for each $j \in \{1,\ldots,m\}$.  There are eight cases to consider.

\vspace{0.2cm}
\textbf{If \boldmath $V_{j,0} > V_{j,1}$.}\,  The optimal strategy segment is \textit{maximum acceleration} to speed $V_{j,0}$, \textit{speedhold} at speed $V_{j,0}$, \textit{maximum acceleration} to speed $v_j(t_1) = W_{j,1} > V_{j,0}$ and \textit{coast} from speed $W_{j,1}$ to speed $V_{j,1}$.   We require $t_1 > t_{a,j}(0,V_{j,0}) + t_{a,j}(V_{j,0},W_{j,1}) = t_{a,j}(0,W_{j,1})$ and $t_2 > t_1 + t_{c,j}(W_{j,1},V_{j,1})$.  The speed $W_{j,1}$ is determined by continuity of $\eta_j$.

\vspace{0.2cm}
\textbf{If \boldmath $V_{j,0} < V_{j,1}$.}\, The optimal strategy segment is \textit{maximum acceleration} to speed $V_{j,0}$, \textit{speedhold} at speed $V_{j,0}$, \textit{coast} to speed $v_j(t_1) = W_{j,1} < V_{j,0}$ and \textit{maximum acceleration} from speed $W_{j,1}$ to speed $V_{j,1}$.   We require $t_1 > t_{a,j}(0,V_{j,0}) + t_{c,j}(V_{j,0},W_{j,1})$ and $t_2 > t_1 + t_{a,j}(W_{j,1},V_{j,1})$.  The speed $W_{j,1}$ is determined by continuity of $\eta_j$.

\vspace{0.2cm}
\textbf{If \boldmath $V_{j,k-1} > V_{j,k} > V_{j,k+1}$.}\,  The optimal strategy segment is \textit{maximum acceleration} from speed $V_{j,k-1}$ to speed $v_j(t_k) = W_{j,k} > V_{j,k-1}$, \textit{coast} from speed $W_{j,k}$ to speed $V_{j,k}$, \textit{speedhold} at speed $V_{j,k}$, \textit{maximum acceleration} from speed $V_{j,k}$ to speed $v_j(t_{k+1}) = W_{j,k+1} > V_{j,k}$ and \textit{coast} from speed $W_{j,k+1}$ to speed $V_{j,k+1}$.  We require $t_k - t_{k-1} > t_{a,j}(V_{j,k-1},W_{j,k})$, $t_{k+1} - t_k > t_{c,j}(W_{j,k}, V_{j,k}) + t_{a,j}(V_{j,k},W_{j,k+1})$ and $t_{k+2} - t_{k+1} > t_{c,j}(W_{j,k+1}, V_{j,k+1})$.  The speeds $W_{j,k}, W_{j,k+1}$ are determined by continuity of $\eta_j$.

\vspace{0.2cm}
\textbf{If \boldmath $V_{j, k-1} > V_{j, k} < V_{j, k+1}$.}\,  The optimal strategy segment is \textit{maximum acceleration} from speed $V_{j,k-1}$ to speed $v_j(t_k) = W_{j,k} > V_{j,k-1}$, \textit{coast} from speed $W_{j,k}$ to speed $V_{j,k}$, \textit{speedhold} at speed $V_{j,k}$, \textit{coast} from speed $V_{j,k}$ to speed $v_j(t_{k+1}) = W_{j,k+1} < V_{j,k}$ and \textit{maximum acceleration} from speed $W_{j,k+1}$ to speed $V_{j, k+1}$.  We require $t_k - t_{k-1} > t_{a,j}(V_{j,k-1},W_{j,k})$, $t_{k+1} - t_k > t_{c,j}(W_{j,k},V_{j,k}) + t_{c,j}(V_{j,k},W_{j,k+1}) = t_{c,j}(W_{j,k}, W_{j,k+1})$ and $t_{k+2} - t_{k+1} > t_{a,j}(W_{j,k+1}, V_{j,k+1})$.  The speeds $W_{j,k}, W_{j,k+1}$ are determined by continuity of $\eta_j$.

\vspace{0.2cm}
\textbf{If \boldmath $V_{j, k-1} < V_{j, k} > V_{j,k+1}$.}\,  The optimal strategy segment is \textit{coast} from speed $V_{j,k-1}$ to speed $v_j(t_k) = W_{j,k} < V_{j,k-1}$, \textit{maximum acceleration} from speed $W_{j,k}$ to speed $V_{j,k}$, \textit{speedhold} at speed $V_{j,k}$, \textit{maximum acceleration} from speed $V_{j,k}$ to speed $v_j(t_{k+1}) = W_{j,k+1} > V_{j,k}$ and \textit{coast} from speed $W_{j,k+1}$ to speed $V_{j,k+1}$.  We require $t_k - t_{k-1} > t_{c,j}(V_{j,k-1},W_{j,k})$, $t_{k+1} - t_k > t_{a,j}(W_{j,k},V_{j,k}) + t_{a,j}(V_{j,k},W_{j,k+1}) = t_{a,j}(W_{j,k},W_{j,k+1})$ and $t_{k+2} - t_{k+1} > t_{c,j}(W_{j,k+1}, V_{j,k+1})$.  The speeds $W_{j,k}, W_{j, k+1}$ are determined by continuity of $\eta_j$.

\vspace{0.2cm}
\textbf{If \boldmath $V_{j, k-1} < V_{j, k} < V_{j, k+1}$.}\,  The optimal strategy segment is \textit{coast} from speed $V_{j,k-1}$ to speed $v_j(t_k) = W_{j,k} < V_{j,k-1}$, \textit{maximum acceleration} from speed $W_{j,k}$ to speed $V_{j,k}$, \textit{speedhold} at speed $V_{j,k}$, \textit{coast} from speed $V_{j,k}$ to speed $v_j(t_{k+1}) = W_{j,k+1} < V_{j,k}$ and \textit{maximum acceleration} from speed $W_{j,k+1}$ to speed $V_{j,k+1}$.  We require $t_k - t_{k-1} > t_{c,j}(V_{j,k-1},W_{j,k})$, $t_{k+1} - t_k > t_{a,j}(W_{j,k}, V_{j,k}) + t_{c,j}(V_{j,k},W_{j,k+1})$ and $t_{k+2} - t_{k+1} > t_{a,j}(W_{j,k+1}, V_{j,k+1})$.  The speeds $W_{j,k}, W_{j,k+1}$ are determined by continuity of $\eta_j$.

\vspace{0.2cm}
\textbf{If \boldmath $V_{j, n-1} < V_{j, n}$.}\,  The final optimal strategy segment is \textit{coast} from speed $V_{j,n-1}$ to speed $v_j(t_n) = W_{j,n} < V_{j,n-1}$, \textit{maximum acceleration} from speed $W_{j,n}$ to speed $V_{j,n}$, \textit{speedhold} at speed $V_{j,n}$, \textit{coast} from speed $V_{j,n}$ to the optimal braking speed $U_{j,n} = \psi(V_{j,n})/ \varphi^{\, \prime}(V_{j,n}) = V_{j,n} - \varphi(V_{j,n})/\varphi^{\, \prime}(V_{j,n})$ and \textit{maximum brake}.  We require $t_n - t_{n-1} > t_{c,j}(V_{j,n-1},W_{j,n})$ and $t_{n+1} - t_n > t_{a,j}(W_{j,n},V_{j,n}) + t_{c,j}(V_{j,n},U_{j,n}) + t_{b,j}(U_{j,n}, 0)$.

\vspace{0.2cm}
\textbf{If \boldmath $V_{j, n-1} > V_{j,n}$.}\,  The final optimal strategy segment is \textit{maximum acceleration} from speed $V_{j,n-1}$ to speed $v_j(t_n) = W_{j,n} > V_{j,n-1}$, \textit{coast} from speed $W_{j,n}$ to speed $V_{j,n}$, \textit{speedhold} at speed $V_{j,n}$, \textit{coast} from speed $V_{j,n}$ to the optimal braking speed $U_{j,n} = \psi(V_{j,n})/ \varphi^{\, \prime}(V_{j,n}) = V_{j,n} - \varphi(V_{j,n})/\varphi^{\, \prime}(V_{j,n})$ and \textit{maximum brake}.  We require  $t_n - t_{n-1} > t_{a,j}(V_{j,n-1},W_{j,n-1})$ and $t_{n+1} - t_n > t_{c,j}(W_{j,n},V_{j,n}) + t_{c,j}(V_{j,n},U_{j,n}) + t_{b,j}(U_{j,n}, 0) = t_{c,j}(W_{j,n},U_{j,n}) + t_{bj(}U_{j,n},0)$.

\subsection{The transition speeds for \boldmath ${\mathfrak T}_j$ on level track}
\label{ss:tslt}

The speed $v_j(t_k) = W_{j,k}$ for train ${\mathfrak T}_j$ depends only on the optimal driving speeds $V_{j,k-1}$ for $t \in (t_{k-1},t_k)$ and $V_{j,k}$ for $t \in (t_k, t_{k+1})$.  We wish to investigate this dependence.  We will assume that $H_a(v) = A/v$ for some constant $A > 0$.  There are two cases to consider.

\vspace{0.2cm}
\textbf{If \boldmath $W_{j,k} > V_{j,k-1} > V_{j,k}$.}\  In this case continuity of the adjoint variable $\eta_j$ at $(v, t) = (W_{j,k},t_k)$ implies
$$
\frac{\varphi^{\, \prime}(V_j)}{\varphi^{\, \prime}(V_{j,k-1})} \cdot \frac{A - L_{V_{j,k-1}}(W_{j,k})}{A - \varphi(W_{j,k})} = \frac{\varphi^{\, \prime}(V_j)}{\varphi^{\, \prime}(V_{j,k})} \cdot \frac{L_{V_{j,k}}(W_{j,k})}{\varphi(W_{j,k})}
$$
which we can rearrange in the equivalent form
\begin{eqnarray}
\label{etacon1}
\lefteqn{ \varphi^{\, \prime}(V_{j,k}) \varphi(W_{j,k}) [A - L_{V_{j,k-1}}(W_{j,k})] } \hspace{1cm} \nonumber \\
& = & \varphi^{\, \prime}(V_{j,k-1})[A - \varphi(W_{j,k})] L_{V_{j,k}}(W_{j,k}).
\end{eqnarray}
Differentiation of (\ref{etacon1}) with respect to $V_{j, k-1}$ shows that
$$
\frac{\partial W_{j,k}}{\partial V_{j,k-1}} = \frac{N_1(V_{k-1},V_k,W_k)}{D(V_{k-1},V_k, W_k)}
$$
where we have defined
\begin{eqnarray*}
N_1(v_1,v_2,w) & = & \varphi^{\, \prime \prime}(v_1)\left\{ \rule{0cm}{0.35cm} [A - \varphi(w)] L_{v_2}(w) \right. \\
& & \hspace{1cm}  \left. \rule{0cm}{0.35cm} + \varphi(w) \varphi^{\, \prime}(v_2)(w - v_1) \right\} \\
& > & 0
\end{eqnarray*}
and
\begin{eqnarray*}
D(v_1,v_2,w) & = & \varphi^{\, \prime}(v_2) \left\{ \rule{0cm}{0.35cm} \varphi^{\, \prime}(w)[A - L_{v_1}(w)] \right. \\
& & \hspace{2cm} \left. \rule{0cm}{0.35cm} - \varphi^{\, \prime}(v_1)[A - \varphi(w)] \right\} \\
& & + \varphi^{\, \prime}(v_1) \left\{ \rule{0cm}{0.35cm} \varphi^{\, \prime}(w)L_{v_2}(w) - \varphi^{\, \prime}(v_2) \varphi(w) \right\} \\
& > & 0
\end{eqnarray*}
when $w > v_1 > v_2$.  To show that $D > 0$ we argue as follows.  We have
$$
\varphi^{\, \prime}(w)[A - L_{v_1}(w)] - \varphi^{\, \prime}(v_1)[A - \varphi(w)] > 0
$$
because $\varphi^{\, \prime}(w) > \varphi^{\, \prime}(v_1)$ and $A - L_{v_1}(w) > A - \varphi(w)$ when $w > v_1$.   If we define
$$
f_1(w) = \varphi^{\, \prime}(w)L_{v_2}(w) - \varphi^{\, \prime}(v_2) \varphi(w)
$$
then $f_1(v_2) = 0$ and since $f_1^{\, \prime}(w) = \varphi^{\, \prime \prime}(w)L_{v_2}(w) > 0$ it follows that $f_1(w) > 0$ when $w > v_2$.   Therefore $\varphi^{\, \prime}(w)L_{v_2}(w) - \varphi^{\, \prime}(v_2) \varphi(w) > 0$ and hence $D(v_1,v_2,w) > 0$. Consequently
\begin{equation}
\label{dwdv1}
\frac{\partial W_{j,k}}{\partial V_{j,k-1}} > 0.
\end{equation}
Differentiation of (\ref{etacon1}) with respect to $V_k$ shows that
$$
\frac{\partial W_{j,k}}{\partial V_{j,k}} = \frac{N_2(V_{k-1},V_k,W_k)}{D(V_{k-1},V_k, W_k)}
$$
where we have defined
\begin{eqnarray*}
N_2(v_1,v_2,w) & = & \varphi^{\, \prime \prime}(v_2)\left\{ \rule{0cm}{0.35cm} - \varphi(w) [A - L_{v_1}(w)] \right. \\
& & \hspace{1cm}  \left. \rule{0cm}{0.35cm} + [A - \varphi(w)] \varphi^{\, \prime}(v_1)(w - v_2) \right\} \\
& < & - \varphi^{\, \prime \prime}(v_2)[ A - \varphi(w)] \cdot \\
& & \hspace{1cm}  \left\{ \rule{0cm}{0.35cm} \varphi(w) + \varphi^{\, \prime}(v_1)(v_2 - w) \right\} \\
& < & - \varphi^{\, \prime \prime}(v_2)[ A - \varphi(w)] \cdot \\
& & \hspace{1cm}  \left\{ \rule{0cm}{0.35cm} \varphi(w) + \varphi^{\, \prime}(w)(v_2 - w) \right\} \\
& = & - \varphi^{\, \prime \prime}(v_2)[ A - \varphi(w)] L_w(v_2) \\
& < & 0
\end{eqnarray*}
when $w > v_1 > v_2$.  We already know that $D(v_1,v_2,W) > 0$ and so it follows that
\begin{equation}
\label{dwdv2}
\frac{\partial W_{j,k}}{\partial V_{j,k}} < 0.
\end{equation}

\begin{figure}[ht]
\begin{center}
\includegraphics[width=0.95\textwidth]{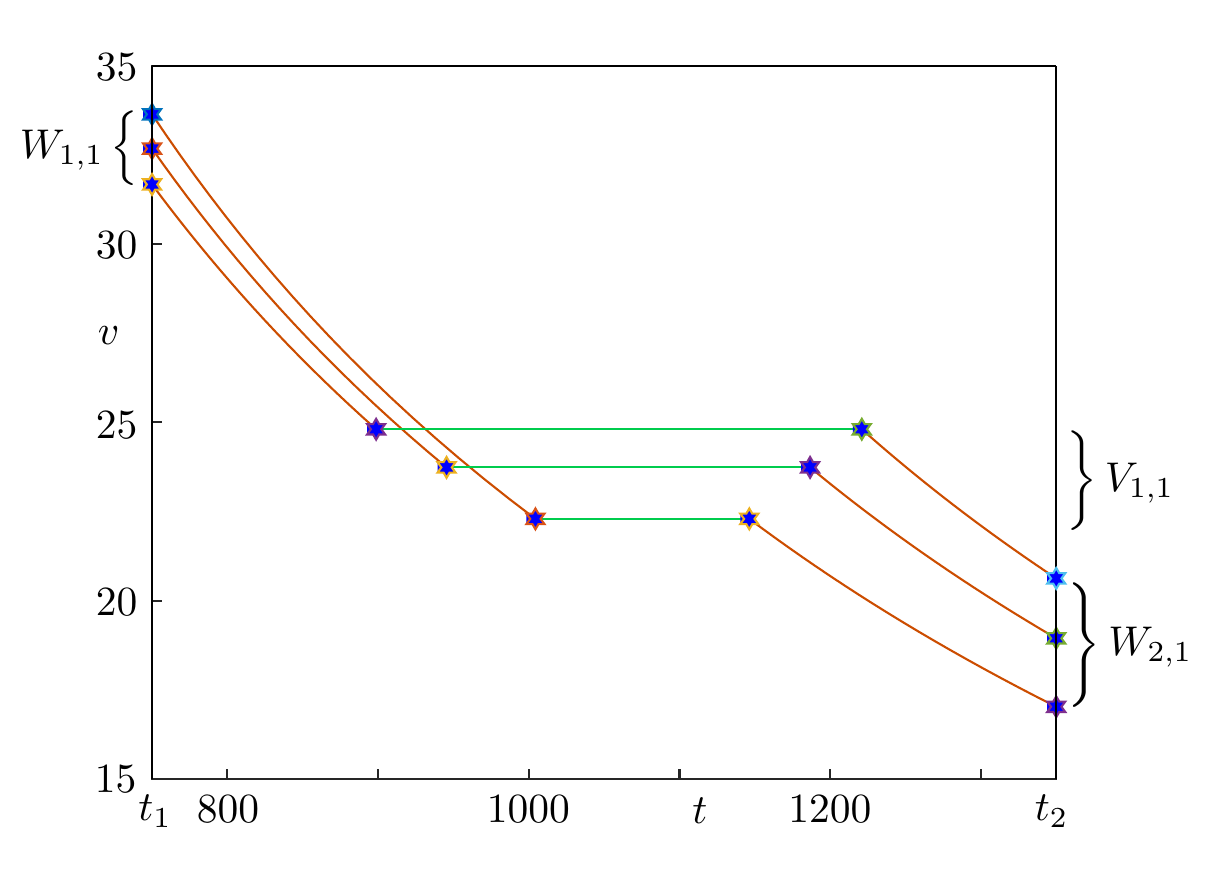}
\end{center}
\vspace{-0.5cm}
\caption{Optimal speed profiles for train ${\mathfrak T}_1$ when $t \in [t_1,t_2] = [750, 1350]$ in Example~\ref{ex:1} with ${\mathcal Q}_1 = 100 M, 200 M, 300 M$. When the restricted optimal driving speed $V_{1,1}$ increases the speed $W_{1,1} = v_1(t_1)$ decreases ($\partial W_{1,1}/\partial V_{1,1} < 0$) and the speed  $W_{1,2} = v_1(t_2)$ increases ($\partial W_{1,2}/\partial V_{1,1} > 0$).  The duration $\Delta \tau_{\,1,1}$ of the speedhold phase also increases.  In this case $E_1 = M \Delta \tau_{\,1,1} \varphi(V_{1,1}) = {\mathcal Q}_1$.}
\label{fig2}
\end{figure}

\begin{figure}[ht]
\begin{center}
\includegraphics[width=0.95\textwidth]{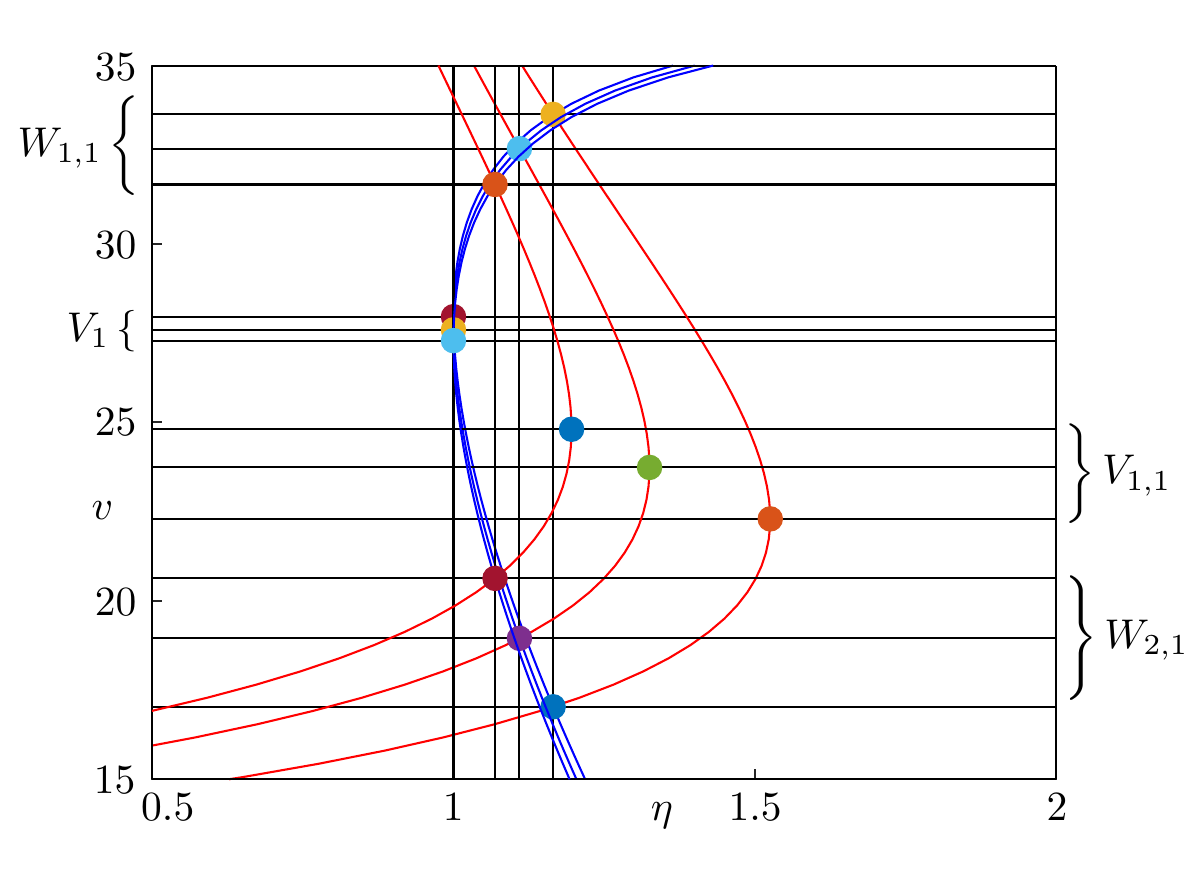}
\end{center}
\vspace{-0.5cm}
\caption{Graphs of the optimal speeds $v_1$ against the modified adjoint variable $\eta_1$ for train ${\mathfrak T}_1$ when $t \in [t_1,t_2] = [750, 1350]$ in Example~\ref{ex:1} with ${\mathcal Q}_1 = 100 M, 200 M, 300 M$ showing (i) the optimal unrestricted driving speeds $V_1$, (ii) the optimal restricted driving speeds $V_{1,1}$, (iii) the switching speeds $W_{1,1} = v_1(t_1)$ and $W_{1,2} = v_1(t_2)$, (iv) the $(\eta_{1,a},v_1)$ curves for the \textit{maximum acceleration} phases on $(0,t_1)$ and $(t_2,T)$ and the associated minimum turning points for $\eta_{1,a}$ at $(1, V_1)$ and (v) the $(\eta_{1,c},v_1)$ curves for the \textit{coast} phases on $(t_1,t_2)$ and the associated maximum turning points for $\eta_{1,c}$ at $(1+w_1, V_{1,1})$.} 
\label{fig3}
\end{figure}

\vspace{0.2cm}
\textbf{If \boldmath $W_{j,k} < V_{j,k-1} < V_{j,k}$.}\,  We can use similar arguments to show that both inequalities (\ref{dwdv1}) and (\ref{dwdv2}) remain true. 

\subsection{Energy consumption for \boldmath ${\mathfrak T}_j$ when $t \in (t_k, t_{k+1})$ on level track}
\label{ss:eclt}

We can write down explicit formul{\ae} for the energy consumed by train ${\mathfrak T}_j$ on the interval $t \in (t_k, t_{k+1})$.  There are eight cases to consider.  In general we can show that $\partial E_{j,k}/\partial V_{j,k} > 0$ and $\partial E_{j,k}/\partial V_{j,k \pm1} < 0$.  We will calculate these derivatives in selected cases to verify the above assertion. The other cases are left to the reader.

\vspace{0.2cm}
\textbf{If \boldmath $V_{j,0} > V_{j,1}$.}\,  The energy consumed on $t \in (t_0, t_1)$ is given by
\begin{eqnarray*}
E_{j,0} & = & \mbox{$\int_0^{W_{j,1}}$} Avdv/[A - \varphi(v)] \\
& & + \left\{ t_1 - t_0 - \mbox{$\int_0^{W_{j,1}}$} vdv/[A - \varphi(v)] \right\} \varphi(V_{j,0}).
\end{eqnarray*}
We know that $W_{j,1}$ depends only on $V_{j,0}$ and $V_{j,1}$.  It follows that $E_{j,0}$ also depends only on $V_{j,0}$ and $V_{j,1}$.

\vspace{0.2cm}
\textbf{If \boldmath $V_{j,0} < V_{j,1}$.}\,  The energy consumed on $t \in (t_0, t_1)$ is given by
\begin{eqnarray*}
E_{j,0} & = & \mbox{$\int_0^{V_{j,0}}$} Avdv/[A - \varphi(v)] \\
& & + \left\{ t_1 - t_0 - \mbox{$\int_0^{W_{j,1}}$} vdv/[A - \varphi(v)] \right.  \\
& & \hspace{1.65cm} \left. - \mbox{$\int_{W_{j,1}}^{V_{j,0}}$} vdv/ \varphi(v) \right\} \varphi(V_{j,0}).
\end{eqnarray*}
Once again $E_{j,0}$ depends only on $V_{j,0}$ and $V_{j,1}$.

\vspace{0.2cm}
\textbf{If \boldmath $V_{j,k-1} > V_{j,k} > V_{j,k+1}$.}\,  We have $W_{j,k} > V_{j,k} < W_{j,k+1}$.  The energy consumed when $t \in (t_k, t_{k+1})$ is
\begin{eqnarray*}
E_{j,k} & = & \mbox{$\int_{V_{j,k}}^{W_{j,k+1}}$} Avdv/[A - \varphi(v)] \\
& & + \left\{ t_{k+1} - t_k - \mbox{$\int_{V_{j,k}}^{W_{j,k}}$} vdv/\varphi(v) \right.  \\
& & \hspace{1.65cm} \left. - \mbox{$\int_{V_{j,k}}^{W_{j,k+1}}$} vdv/ [A - \varphi(v)] \right\} \varphi(V_{j,k}).
\end{eqnarray*}
In this case $E_{j,k}$ depends only on $V_{j,k-1}, V_{j,k}, V_{j,k+1}$ and we calculate
$$
\frac{\partial E_{j,k}}{\partial V_{j,k-1}} = (-1) \frac{W_{j,k}}{\varphi(W_{j,k})} \frac{\partial W_{j,k}}{\partial V_{j,k-1}} < 0,
$$
because $\partial W_{j,k}/\partial V_{j,k-1} > 0$,
\begin{eqnarray*}
\frac{\partial E_{j,k}}{\partial V_{j,k}} & = & \frac{W_{j,k+1} [A - \varphi(V_{j,k})]} {[A - \varphi(W_{j,k+1})]} \frac{\partial W_{j,k+1}}{\partial V_{j,k}} \\
& & \hspace{1cm} - \frac{W_{j,k} \varphi(V_{j,k}}{\varphi(W_{j,k})} \frac{\partial W_{j,k}}{\partial V_{j,k}} \\
& & + \left\{ t_{k+1} - t_k - \mbox{$\int_{V_{j,k}}^{W_{j,k}}$} vdv/\varphi(v) \right.  \\
& & \hspace{1.65cm} \left. - \mbox{$\int_{V_{j,k}}^{W_{j,k+1}}$} vdv/ [A - \varphi(v)] \right\} \varphi^{\, \prime}(V_{j,k}) \\
& > & 0
\end{eqnarray*}
because $\partial W_{j,k+1}/\partial V_{j,k} > 0$ and $\partial W_{j,k}/\partial V_{j,k} < 0$, and
$$
\frac{\partial E_{j,k}}{\partial V_{j,k+1}} = \frac{W_{j,k+1}[A - \varphi(V_{j,k}]}{[A - \varphi(W_{j,k+1}]} \frac{\partial W_{j,k+1}}{\partial V_{j,k+1}} < 0
$$
because $\partial W_{j,k+1}/\partial V_{j,k+1} < 0$.

\vspace{0.2cm}
\textbf{If \boldmath $V_{j,k-1} > V_{j,k} < V_{j,k+1}$.}  We have $W_{j,k} > V_{j,k} > W_{j,k+1}$.  The energy consumed when $t \in (t_k, t_{k+1})$ is
$$
E_{j,k} = \left\{ t_{k+1} - t_k - \mbox{$\int_{W_{j,k+1}}^{W_{j,k}}$} vdv/\varphi(v) \right\} \varphi(V_{j,k}).
$$
We can see that $E_{j,k}$ depends only on $V_{j,k-1}, V_{j,k}, V_{j,k+1}$.  In this case we calculate
$$
\frac{\partial E_{j,k}}{\partial V_{j,k-1}} = (-1) \frac{W_{j,k}}{\varphi(w_{j,k})} \frac{\partial W_{j,k}}{\partial V_{j,k-1}} < 0
$$
because $\partial W_{j,k}/\partial V_{j,k-1} > 0$,
$$
\frac{\partial E_{j,k}}{\partial V_{j,k+1}} = \left\{ t_{k+1} - t_k - \mbox{$\int_{W_{j,k+1}}^{W_{j,k}}$} vdv/\varphi(v) \right\} \varphi^{\, \prime}(V_{j,k}) > 0,
$$
and
$$
\frac{\partial E_{j,k}}{\partial V_{j,k+1}} = (-1) \frac{W_{j,k+1}}{\varphi(w_{j,k+1})} \frac{\partial W_{j,k+1}}{\partial V_{j,k+1}} < 0
$$
because $\partial W_{j,k}/\partial V_{j,k+1} > 0$.  This case is illustrated in Figures~\ref{fig2} and \ref{fig3}.

\vspace{0.2cm}
\textbf{If \boldmath $V_{j,k-1} < V_{j,k} > V_{j,k+1}$.}\,  The energy consumed when $t \in (t_k, t_{k+1})$ is given by
\begin{eqnarray*}
E_{j,k} & = & \mbox{$\int_{W_{j,k}}^{V_{j,k}}$} Avdv/[A - \varphi(v)] \\
& & + \left\{ t_{k+1} - t_k - \mbox{$\int_{W_{j,k}}^{V_{j,k}}$} vdv/[A - \varphi(v)] \right.  \\
& & \hspace{1.65cm} \left. - \mbox{$\int_{W_{j,k+1}}^{V_{j,k}}$} vdv/ \varphi(v) \right\} \varphi(V_{j,k}).
\end{eqnarray*}
Therefore $E_{j,k}$ depends only on $V_{j,k-1}, V_{j,k}, V_{j,k+1}$.

\vspace{0.2cm}
\textbf{If \boldmath $V_{j,k-1} < V_{j,k} < V_{j,k+1}$.}\,  The energy consumed when $t \in (t_k, t_{k+1})$ is given by
\begin{eqnarray*}
E_{j,k} & = & \mbox{$\int_{V_{j,k}}^{W_{j,k+1}}$} Avdv/[A - \varphi(v)] \\
& & + \left\{ t_{k+1} - t_k - \mbox{$\int_{V_{j,k}}^{W_{j,k}}$} vdv/\varphi(v) \right.  \\
& & \hspace{1.65cm} \left. - \mbox{$\int_{V_{j,k}}^{W_{j,k+1}}$} vdv/ [A - \varphi(v)] \right\} \varphi(V_{j,k}).
\end{eqnarray*}
Once again $E_{j,k}$ depends only on $V_{j,k-1}, V_{j,k}, V_{j,k+1}$.

\vspace{0.2cm}
\textbf{If \boldmath $V_{j,n-1} < V_{j,n}$.}\,  The energy consumed when $t \in (t_n, t_{n+1})$ is given by
\begin{eqnarray*}
E_{j,n} & = & \mbox{$\int_{W_{j,n}}^{V_{j,n}}$} Avdv/[A - \varphi(v)] \\
& & \hspace{-1.2cm} + \left\{ t_{n+1} - t_n - \mbox{$\int_{W_{j,n}}^{V_{j,n}}$} vdv/[A - \varphi(v)] \right.  \\
& & \hspace{-1.2cm}  \left. - \mbox{$\int_{U_{j,n}}^{V_{j,n}}$} vdv/\varphi(v) - \mbox{$\int_0^{U_{j,n}}$} vdv/[H_b(v) + \varphi(v)] \right\} \varphi(V_{j,n}).
\end{eqnarray*}
Clearly $E_{j,n}$ depends only on $V_{j,n-1}$ and $V_{j,n}$.

\vspace{0.2cm}
\textbf{If \boldmath $V_{j,n-1} > V_{j,n}$.}\,  The energy consumed when $t \in (t_n, t_{n+1})$ is given by
\begin{eqnarray*}
E_{j,n} & = & \left\{ t_{n+1} - t_n - \mbox{$\int_{U_{j,n}}^{W_{j,n}}$} vdv/\varphi(v) \right.\\
& & \hspace{1.5cm} \left. - \mbox{$\int_0^{U_{j,n}}$} vdv/[H_b(v) + \varphi(v)] \right\} \varphi(V_{j,n}).
\end{eqnarray*}
Thus $E_{j,n}$ depends only on $V_{j,n-1}$ and $V_{j,n}$.

\subsection{A feasible optimal strategy}
\label{ss:fop}

For each fixed value of $\bfw \geq \bfzero$ the optimal strategy for train ${\mathfrak T}_j$ takes precisely the same form as the optimal strategy obtained by Howlett \textit{et al.}~\cite{how14} for a single train ${\mathfrak T}_1$ with energy consumption constraints on the intervals $(t_k, t_{k+1})$ for each $k \in \{0,\ldots,n\}$.  Consider the problem with a single train ${\mathfrak T}_1$.  The problem could conceivably be solved by the following naive iteration. Suppose a value of $\bfw = (w_1,\ldots,w_n) \geq \bfzero$ is given.  For each choice of the unconstrained optimal driving speed $V_1$ the optimal constrained driving speeds $V_{1,k}$ on the intervals $(t_k, t_{k+1})$ are defined by the formula
\begin{equation}
\label{wkv1v1k}
1 + w_k = \varphi^{\, \prime}(V_1)/\varphi^{\, \prime}(V_{1,k}).
\end{equation}
If $w_{\ell} = 0$ then $\ell \notin {\mathcal A}$ and $V_{1,\ell} = V_1$.  If $w_k > 0$ then $k \in {\mathcal A}$ and $V_{1,k} < V_1$.  Thus an entire strategy of optimal type is determined.   In order to find a feasible strategy of optimal type we must choose $V_1$ so that the distance constraint
\begin{equation}
\label{xv1w}
X = \mbox{$\int_0^T$} v_1(t) dt
\end{equation}
is satisfied.  However we also need to satisfy the individual energy consumption constraints
\begin{equation}
\label{q1kv1wk}
\mbox{$\int_{t_k}^{t_{k+1}}$} u_{a,1}(t)v_1(t) dt - {\mathcal Q}_{1,k} = 0
\end{equation}
on each interval $(t_k,t_{k+1})$ with $k \in {\mathcal A}$.  Consequently the next step is to find a new value $w_k > 0$    for each $k \in {\mathcal A}$ so that the corresponding constraints (\ref{q1kv1wk}) are all satisfied.  This step can be implemented with a multi-dimensional Newton iteration.  The new value of $\bfw \geq \bfzero$ will change the value of $V_{1,k}$ for each $k \in {\mathcal A}$.  Thus it is now necessary to choose a new value of $V_1$ so that the distance constraint (\ref{xv1w}) is satisfied.  This step could be implemented with an elementary Newton iteration.  The new value of $V_1$ generates new values of $V_{1,k}$ using (\ref{wkv1v1k}).  The process must necessarily proceed by alternate adjustment of $\bfw$ and $V_1$ until the constraints (\ref{q1kv1wk}) and (\ref{xv1w}) are all satisfied to sufficient accuracy.  In practice Howlett \textit{et al.}~\cite{how14} suggested that either a multi-dimensional Newton iteration or a modified midpoint method could be used.  Both methods converge rapidly.

In the case of a fleet of trains $\{ {\mathfrak T}_j \}_{j=1}^m$ the naive iteration proceeds in much the same way.  Suppose we choose an initial value of $\bfw \geq \bfzero$.  For each fixed vector $\bfV = (V_1,\ldots,V_m) > \bfzero$ of unconstrained optimal driving speeds the values of the constrained optimal driving speeds $V_{j,k}$ for each $j \in \{1,\ldots,m\}$ and each $k \in {\mathcal A}$ are determined by the formula
\begin{equation}
\label{wkvjvjk}
1 + w_k = \varphi^{\, \prime}(V_j)/\varphi^{\, \prime}(V_{j,k}).
\end{equation}
Thus, for each choice of $\bfV$ an entire strategy of optimal type is determined for each train.   In order to find a feasible strategy of optimal type we must choose $\bfV$ so that the distance constraint
\begin{equation}
\label{xjvjw}
X_j = \mbox{$\int_0^T$} v_j(t) dt
\end{equation}
is satisfied for each $j=1,2,\ldots,m$.  However we also need to satisfy the individual energy consumption constraints
\begin{equation}
\label{qkvwk}
\mbox{$\sum_{j=1}^m \int_{t_k}^{t_{k+1}}$} u_{a,j}(t)v_j(t) dt - {\mathcal Q}_k = 0
\end{equation}
on each interval $(t_k,t_{k+1})$ with $k \in {\mathcal A}$.  Thus we must now adjust the value of $\bfw \geq \bfzero$.  The iteration then proceeds by alternate adjustment of the vectors $\bfV$ and $\bfw$.  In practice one can use a rapidly convergent multi-dimensional Newton iteration or a modified midpoint method.

\begin{figure}[ht]
\begin{center}
\includegraphics[width=0.95\textwidth]{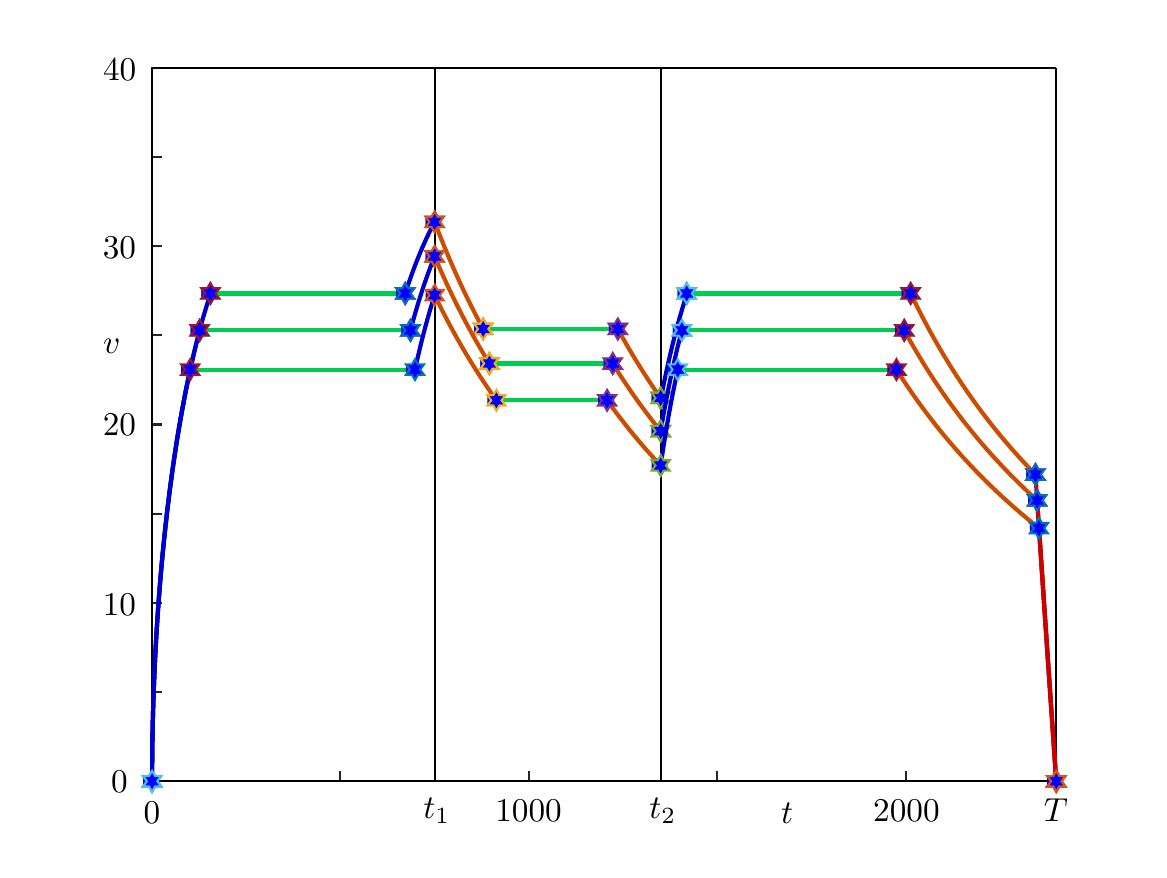}
\end{center}
\vspace{-0.5cm}
\caption{Optimal speed profiles for trains $\{ {\mathfrak T}_j\}_{j=1}^3$ in Example~\ref{ex:2}.} 
\label{fig4}
\end{figure}

\begin{example}
\label{ex:2}
We consider a fleet of three identical trains $\{{\mathfrak T}_j\}_{j=1}^3$ and use the same model as in Example~\ref{ex:1} with $r(v) = r_0 + r_2v^2$ where $r_0 = 6.75 \times 10^{-3}$ and $r_2 = 5 \times 10^{-5}$,  $H_a(v) = A/v$ where $A = 3$ and $H_b(v) = 3 \times 10^{-1}$.  The journey distances are $X_1 = 60000, X_2 = 55000, X_3 = 50000$ m and the journey time is $T= 2400$ s.  We find optimal speed profiles $v_j = v_j(t)$ for each $j=1,2,3$ subject to the single energy consumption constraint $\sum_{j=1}^3 E_{j,1} \leq {\mathcal Q}_1 = (800 M)$~J  on the $10$\textendash minute interval $t \in (t_1,t_2) = (750, 1350)$ s.  The cost of the journey for train ${\mathfrak T}_j$ is the total energy consumption $J_j = M \int_0^T u_{j,a}(t)v_j1(t) dt$ ~J.  The values of the key variables are shown in TABLE~2 and the corresponding optimal speed profiles are shown in Figure~\ref{fig4}.  The optimal weight factor is $w_1 \approx 0.152612$.  This problem was solved by the naive method using alternate \textit{ad hoc} adjustments to the parameters $w_1$ and $\bfV = (V_1,V2,V3)$ until the distance constraints $\int_0^T v_j(t) dt = X_j$ for each $j=1,2,3$ and the intermediate energy consumption constraint $\sum_{j=1}^3 M \int_{t_{j-1}}^{t_j} u_{a,j}(t)v_j(t) dt = {\mathcal Q}_1$ were all satisfied to sufficient accuracy. $\hfill \Box$
\end{example}

\begin{table}[ht]
\vspace{0.25cm}
\begin{center}
\caption{Values of the key variables in Example~\ref{ex:2}.}
\begin{tabular}{|c|c|c|c|c|c|c|c|} \hline
\rule{0cm}{0.25cm} \hspace{-2mm} ${\mathfrak T}_j$ \hspace{-2mm} & $J_j$ & $E_{j,1}$ & $V_j$ & $W_{j,1}$ & $V_{j,1}$ & $W_{j,2}$ & $U_j$ \\ \hline
\rule{0cm}{0.25cm} \hspace{-2mm} $j=1$ \hspace{-2mm} & $2600 M$  & $353 M$ & 27.35  &  31.35  &  25.35  &  21.49  &  17.20 \\ \hline
\rule{0cm}{0.25cm} \hspace{-2mm} $j=2$ \hspace{-2mm} & $2127 M$ & $262 M$ & 25.28  &  29.42  &  23.42  & 19.63  &  15.74 \\ \hline
\rule{0cm}{0.25cm} \hspace{-2mm} $j=3$ \hspace{-2mm} & $1688 M$ & $185 M$ & 23.08  & 27.25  & 21.35  & 17.70  & 14.19 \\ \hline
\end{tabular}
\end{center}
\end{table}

\begin{figure}[ht]
\begin{center}
\includegraphics[width=0.95\textwidth]{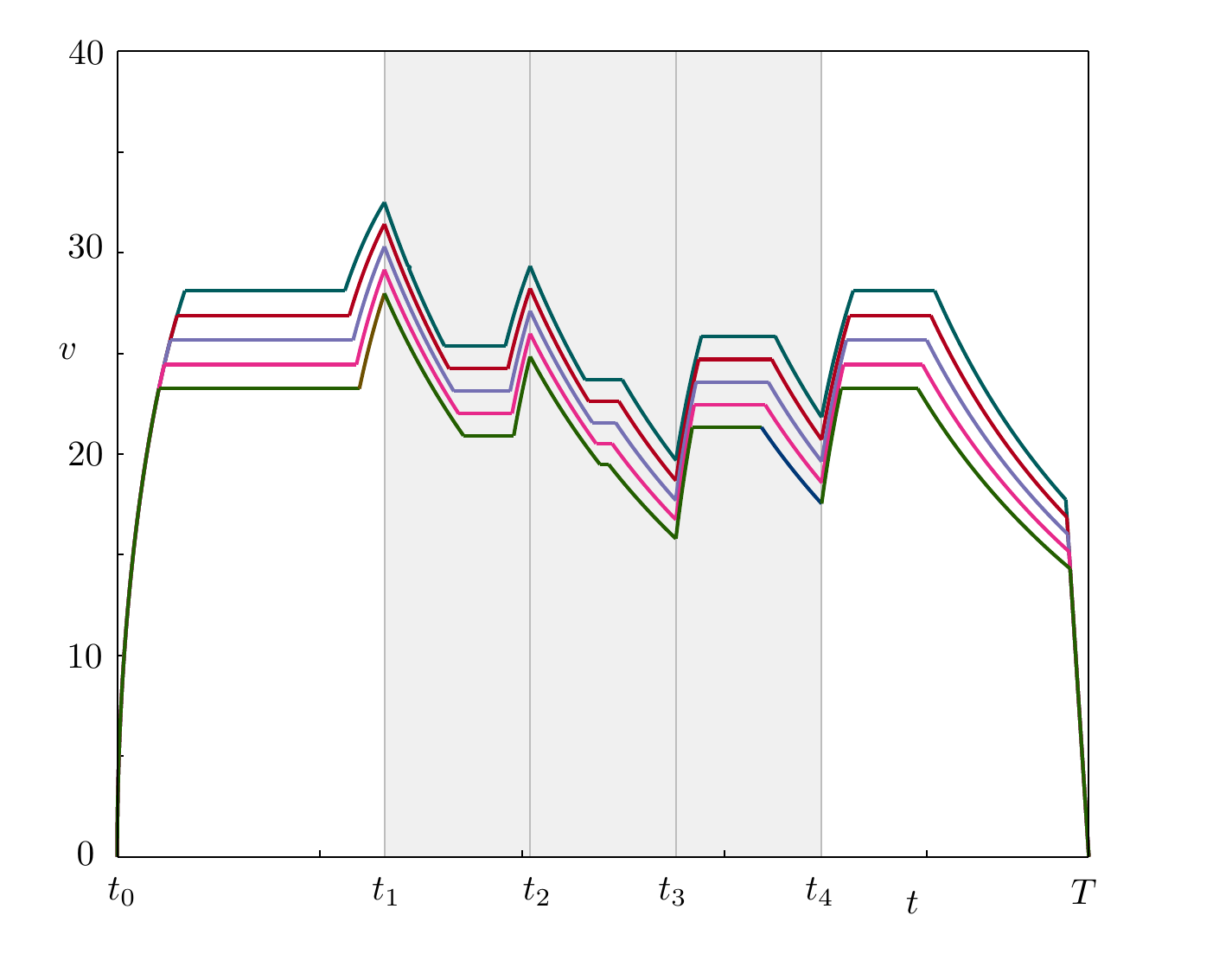}
\end{center}
\vspace{-0.5cm}
\caption{Optimal speed profiles for trains $\{ {\mathfrak T}_j\}_{j=1}^5$ in Example~\ref{ex:3}.} 
\label{fig5}
\end{figure}

\begin{example}
\label{ex:3}
We consider a fleet of five identical trains $\{{\mathfrak T}_j\}_{j=1}^5$ and use the same model as we used in Example~\ref{ex:1}.  We assume each train has mass $M$.  The journey distances are $X_1 = 60000, X_2 = 57500, X_3 = 55000$, $X_4 = 52500$ and $X_5 = 50000$~m and the journey time is $T= 2400$~s.  If there are no intermediate energy constraints the optimal strategy for each train is \textit{maximum acce;leration}, \textit{speedhold} at the optimal driving speed, \textit{coast} to the optimal braking speed and \textit{maximum brake}.  The optimal driving speeds are $\bfV \approx (26.68, 25.54, 24.41, 23.28, 22.16)$, the optimal braking speeds are $\bfU \approx (16.73, 15.93, 15.13, 14.33, 13.54)$ and the optimal journey costs are
$$
\bfJ_{\min} \approx (2541 M, 2268 M, 2018 M, 1787 M, 1577 M)~J.
$$
These strategies are not feasible in the problem we shall consider.  Suppose the fixed intermediate times are
$$
\bft = (0, 660, 1020, 1380, 1740, 2400)~\mbox{$s$}.
$$
We wish to find optimal speed profiles $v_j = v_j(t)$ for each ${\mathfrak T}_j$ subject to the energy constraints
$$
\mbox{$\sum_{j=1}^5$} E_{j,1} \leq {\mathcal Q}_1 = (1300 M)~\mbox{J} 
$$
for $t \in (t_1,t_2)$~s,
$$
\mbox{$\sum_{j=1}^5$} E_{j,2} \leq {\mathcal Q}_2 = (200 M)~\mbox{J}
$$
for $t \in (t_2,t_3)$~s, and
$$
\mbox{$\sum_{j=1}^5$} E_{j,3} \leq {\mathcal Q}_3 = 1500 M~\mbox{J}
$$
for $t \in (t_3,t_4)$~s.  The cost of the journeys for the entire fleet is the total energy consumption
$$
\mbox{$\sum_{j=1}^5$} J_j = M \int_0^T u_{j,a}(t)v_j1(t) dt~\mbox{J}.
$$
The optimal strategies for each train take the general form described in Section~\ref{ss:oslt}.  The values of the key variables are shown in TABLE~3 and the corresponding optimal speed profiles are shown in Figure~\ref{fig5}.  The optimal weight factor is $\bfw \approx (0.213310, 0.378544, 0.170739)$.  The total cost of the optimal strategies with no energy constraints is $\| \bfJ_{\min} \|_1 \approx (10191 M)$ J.  The total cost for the restricted strategies with the given energy constraints is $\| \bfJ \|_1 \approx (10399 M)$ J.  The solution was found using the \textit{fsolve} routine in the high-level, general-purpose programming language \textit{Python}.  The \textit{fsolve} routine is an implementation of \textit{Powell's dog leg method}\textemdash also known as \textit{Powell's hybrid method}\textemdash for solution of multi-dimensional non-linear equations.  The method~\cite{pow1} combines a Gauss\textendash Newton algorithm with gradient descent and uses an explicit trust region. $\hfill \Box$
\end{example}

\begin{table}[ht]
\vspace{0.25cm}
\begin{center}
\caption{Values of the key variables in Example~\ref{ex:3}.}
\begin{tabular}{|c|c|c|c|c|c|} \hline
\rule{0cm}{0.25cm} $  $ & $ {\mathfrak T}_1$ & ${\mathfrak T}_2$ & ${\mathfrak T}_3$ & ${\mathfrak T}_5$ & ${\mathfrak T}_5$ \\ \hline
\rule{0cm}{0.25cm} $V_j$ & 28.11 & 26.88 & 25.65  &  24.44  &  23.2467  \\ \hline
\rule{0cm}{0.25cm} $W_{j,1}$ & 32.51 & 31.42 & 30.31  & 29.16 &  27.98  \\ \hline
\rule{0cm}{0.25cm} $V_{j,1}$ & 25.37 & 24.24 & 23.12  &  22.61  &  20.92  \\ \hline
\rule{0cm}{0.25cm} $W_{j,2}$ & 29.34 & 28.23 & 27.11  & 25.98  & 24.84  \\ \hline
\rule{0cm}{0.25cm} $V_{j,2}$ & 23.69  & 22.62 & 21.56  & 20.52  &  19.48  \\ \hline
\rule{0cm}{0.25cm} $W_{j,3}$ & 19.69 & 18.68 & 17.70  &  16.74  &  15.80  \\ \hline
\rule{0cm}{0.25cm} $V_{j,3}$ & 25.86 & 24.71 & 23.57 & 22.45  & 21.33  \\ \hline
\rule{0cm}{0.25cm} $W_{j,4}$ & 21.83 & 20.71 & 19.63  & 18.57  & 17.54  \\ \hline
\rule{0cm}{0.25cm} $U_j$ & 17.73 & 16.87 & 16.01  & 15.15 & 14.31  \\ \hline
\rule{0cm}{0.25cm} $J_j $ & $2590 M$ & $2314 M$ & $2059 M$ & $1825 M$  & $1611 M$ \\ \hline
\rule{0cm}{0.25cm} $E_{j,1}$ & $332 M$ & $292 M$ & $256 M$  & $224 M$  & $195 M$  \\ \hline
\rule{0cm}{0.25cm} $E_{j,2}$ & $75 M$ & $54 M$ & $37 M$  & $23 M$  & $11 M$  \\ \hline
\rule{0cm}{0.25cm} $E_{j,3}$ & $379 M$ & $335 M$ & $296 M$ & $261 M$  & $229 M$   \\ \hline
\end{tabular}
\end{center}
\end{table}

\section{Conclusions}
\label{s:con}

We have shown that the optimal driving strategy proposed by Howlett \textit{et al.}~\cite{how14} for a single train subject to constraints on energy consumption on fixed intermediate time intervals can be extended to a fleet of trains with each individual driving strategy taking the same general form.  We also showed that our proposed strategies could be calculated precisely in practice for a fleet of trains on level track using standard methods of numerical analysis.  The new algorithms could easily be added to an existing on-board DAS such as \textit{Energymiser}\textsuperscript{\textregistered}.  Our future research will focus on the development of stable algorithms that calculate optimal strategies on tracks with non-zero gradients.  In this regard we note that the Pontryagin analysis in this paper remains valid on tracks with non-zero gradient.


\begin{thebibliography}{99}

\bibitem{alb4} Amie Albrecht, Phil Howlett, Peter Pudney, Xuan Vu (2013), Energy-efficient train control: from local convexity to global optimization and uniqueness, {\em Automatica}, {\bf 49}, 3072-3078. http://dx.doi.org/10.1016/j.automatica.2013.07.008.

\bibitem{alb9} Amie Albrecht, Phil Howlett, Peter Pudney, Xuan Vu, Peng Zhou, 2015, The key principles of optimal train control---Part 1: Formulation of the model, strategies of optimal type, evolutionary lines, location of optimal switching points, {\em Transportation Research Part B: Methodological}, {\bf 94}, 482-508.  DOI: http://dx.doi.org/10.1016/j.trb.2015.07.023.

\bibitem{alb10} Amie Albrecht, Phil Howlett, Peter Pudney, Xuan Vu, Peng Zhou, 2015, The key principles of optimal train control---Part 2: Existence of an optimal strategy, the local energy minimization principle, uniqueness, computational techniques, {\em Transportation Research Part B: Methodological}, {\bf 94}, 509-538.  DOI: http://dx.doi.org/10.1016/j.trb.2015.07.024.

\bibitem{alb13} Amie Albrecht, Phil Howlett, Peter Pudney, (2020). Calculation of Optimal Driving Strategies for Two Successive Trains with Safe Separation, IEEE Transactions on Intelligent Transportation Systems, {\bf 23} (1), 280\textendash 295. https://doi.org/10.1109/TITS.2020.3010245.

\bibitem {asn1} I.A.~Asnis, A.V.~Dmitruk, N.P.~Osmolovskii, 1985, Solution of the problem of the energetically optimal control of the motion of a train by the maximum principle, {\em U.S.S.R. Comput.Maths.Math.Phys.}, {\bf 25}, 6, 37\textendash 44. DOI: https://doi.org/10.1016/0041-5553(85)90006-0.

\bibitem{bar1}L.A.~Baranov, E.~Erofeyev, I.~Golovitcher, V.~Maksimov, {\em Automated control for Electric Locomotives and Multiple Units}. Monograph, Transport: Moscow, 1990.

\bibitem{bar2}L.A.~Baranov, I.S.~Meleshin, LM.~Chin, 2011, Optimal control of a subway train with regard to the criteria of minimum energy consumption, {\em Russian Electrical Engineering},  {\bf 82}, 8, 405\textendash 410.  DOI: https://doi.org/10.3103/S1068371211080049.

\bibitem{bur1}R.~L.~Burdett, E.~Kozan, Techniques for inserting additional trains into existing timetables. Transportation Research B, 43(8), (2009), 821-836.  DOI: https://doi.org/10.1016/j.trb.2009.02.005.

\bibitem{bur2}R.~L.~Burdett, E.~Kozan, A sequencing approach for train timetabling, {\em OR Spectrum}, 32(1), (2010), 163-193.  DOI: 10.1007/s00291-008-0143-6.

\bibitem{cap1} A.~Caprara, M.~Fischetti, P.~Toth, 2002. Modeling and solving the train timetabling problem. {\em Oper. Res.}, {\bf 50} (5), 851\textendash 861.  DOI: https://doi.org/ 10.1287/opre.50.5.851.362.

\bibitem{che1} J.~Cheng, P.G.~Howlett, 1992,  Application of critical velocities to the minimisation of fuel consumption in the control of trains, {\em Automatica}, {\bf 28}, 1, 165\textendash 169.  DOI: https://doi.org/10.1016/0005-1098(92)90017-A.

\bibitem{cho1} Chou, M. and Xia, X., 2007, Optimal cruise control of heavy-haul trains equipped with electronically controlled pneumatic brake systems, {\em Control Engineering Practice}, {\bf 15}, 511--519.

\bibitem{dav1} WJ Davis, Jr. 1926, The tractive resistance of electric locomotives and cars, {\em General Electric Review}, {\bf 29}, 2\textendash24.

\bibitem{dor1} M.J Dorfman and J. Medanic, 2004, Scheduling trains on a railway network using a discrete event model of railway traffic, Transportation Research Part B: Methodological, {\bf 38} (1), 81--98.

\bibitem{dub1}  A.Ya. Dubovitskii and A.A. Milyutin, 1965, Extremum problems in the presence of restrictions, USSR Comput. Math. and Math. Phys., {\bf 5} (3), 1--80.

\bibitem{gir1} I.V. Girsanov, 1972,  Lectures on Mathematical Theory of Extremum Problems,  Lecture Notes in Economics and Mathematical Systems, Eds. M. Beckman, G. Goos and H.P. K\"{u}nzi, {\bf 67}, Springer-Verlag.

\bibitem{gol2}I.~Golovitcher, 1986, control algorithms for automatic operation of rail vehicles. {\em Automated and Remote control. Journal of Russian (USSR) Academy of Science} (Automatika i Telemekhanika) {\bf 11}, 118\textendash 126.
\bibitem{gup1} Shuvomoy Das Gupta, J.~Kevin Tobin, Lacra Pavel, 2016, A two-step linear programming model for energy-efficient timetables in metro railway networks, Transportation Research Part B: Methodological, {\bf 93}, 57\textendash 74. DOI: https://doi.org/10.1016/j.trb.2016.07.003.

\bibitem{har1} Richard F. Hartl, Suresh P. Sethi and Raymond G. Vickson, 1995, A survey of the maximum principles for optimal control problems with state constraints, SIAM Review, {\bf 37} (2), 181--218. 

\bibitem{how1}P.G.~Howlett, 1990, An optimal strategy for the control of a train, {\em J. Aust. Math. Soc., Series B}, (now {\em ANZIAM J.}), {\bf 31}, 454\textendash 471.  DOI:  https:// doi.org/10.1017/ S0334270000006780.

\bibitem{how3}P.G.~Howlett, P.J.~Pudney, I.P.~Milroy, 1994, Energy-efficient train control, {\em control Engineering Practice}, {\bf 2}, 2, 193\textendash 200.  DOI: https://doi.org/10.1016/0967-0661(94)90198-8.

\bibitem{how4} P.G.~Howlett, P.J.~Pudney, {\it Energy-Efficient Train control}, Advances in Industrial control, Springer, London, 1995.

\bibitem{how6} Phil Howlett, Cheng Jiaxing, 1997, Optimal Driving Strategies for a Train on a Track with Continuously Varying Gradient, {\it ANZIAM J.} formerly {\it J. Aust. Math. Soc. Ser. B}, {\bf 38}, 388\textendash410.  DOI: https://doi.org/10.1017/ S0334270000000746.

\bibitem {how7} Phil Howlett, 2000, The optimal control of a train, Ann Oper Res, {\bf 98}, 65--87.

\bibitem{how8} Howlett, P.G., Leizarowitz, A., 2001. Optimal strategies for vehicle control problems with finite control sets. Dynamics Continuous Discrete and Impulsive Systems, Series B {\bf 8}, 41\textendash 69.

\bibitem{how11} Phil Howlett, Peter Pudney, Amie Albrecht (2022). Optimal Driving Strategies for a Fleet of Trains on Level Track with Prescribed Intermediate Signal Times and Safe Separation.  {\em Transportation Science}, {\bf 57} (2),   399\textendash 423. https:// doi.org/10.1287/trsc.2022.1170.

\bibitem{how12} Phil Howlett, Peter Pudney, (2023). The cost differential for an optimal train journey on level track, Journal of Rail Transport, Planning and Management (JRTPM), {\bf 26}, 100393, 1\textendash 21. doi.org/10.1016/ j.jrtpm.2023.100393.

\bibitem{how13} Phil Howlett, Peter Pudney, Amie Albrecht, (2023). An optimal train journey with bounds on energy consumption during specified intermediate time intervals, Journal of Rail Transport, Planning and Management (JRTPM), {\bf 27}, 100391, 1\textendash 18.  doi.org/10.1016/ j.jrtpm.2023.100391.

\bibitem{how14} Phil Howlett, Maria Kapsis, Peter Pudney, (2025).  Optimal driving strategies for trains on level track with bounds on energy consumption on specified intermediate time intervals, Journal of Rail Transport Planning and Management, \textbf{36}, 100550, 18 pages.  https://doi.org/10.1016/j.jrtpm.2025.100550. Preprint available on SSRN at https:// dx.doi.org/10.2139/ssrn.5084990.


\bibitem{kap1} Maria Kapsis, Peter Pudney, Phil Howlett, Peng Zhou, Amie Albrecht, (2025).  Reducing energy during specified time intervals for multiple trains using simple train models, ANZIAM J, \textbf{67}, e30, 16 pages (published online, August 26, 2025), https://doi.org/10.1017/S1446181125100138.  

\bibitem {khm1} Eugene Khmelnitsky, 2000, On an Optimal Control Problem of Train Operation, IEEE T Automat Contr, {\bf 45}, {\bf 7}, 1257--1266.

\bibitem {liu1} Liu, R. and Golovitcher, I., 2003, Energy-efficient operation of rail vehicles, Transport Res A-Pol, {\bf 37}, 917--932.

\bibitem{lilo1}X.~Li, H.K.~Lo, 2014,  An energy-efficient scheduling and speed control approach for metro rail operations. {\em Transportation Research Part B: Methodological}, {\bf 64}, 73\textendash 89.  DOI:  https://doi.org/10.1016/j.trb.2014.03.006.

\bibitem{lilo2} X.~Li, H.K.~Lo, 2014,  Energy minimization in dynamic train scheduling and control for metro rail operations. {\em Transportation Research Part B: Methodological}, {\bf 70}, 269\textendash 284.  DOI: https://doi.org/10.1016/ j.trb.2014.09.009.

\bibitem{lius2} Liu S., Kozan E., 2011, Scheduling trains with priorities: a no-wait blocking parallel-machine job-shop scheduling model, Transportation Science, {\bf 45} (2),175--198.

\bibitem{lue1} D.G. Luenberger, Optimization by Vector Space Methods, John Wiley \& Sons, 1969.

\bibitem{pam1} Abdoulaye Pam, Tony Letrouv\'{e}, Peng Zhou, Robert Yee, Peter Pudney, Olivier Grellier, 2022.   Targeted Traction Power Modulation of High-Speed Trains for Stabilization of the Electric Supply Network with Electric Flexibility,  {\em Proceedings IEEE Vehicle Power and Propulsion Conference (VPPC)}, 1\textendash 4 November, 2022.  doi.org/10.1109/VPPC55846.2022 .10003344.

\bibitem{pow1} MJD Powell (1970). A New Algorithm for Unconstrained Optimization. In: JB Rosen, OL Mangasarian and K Ritter Eds., \textit{Nonlinear Programming}, Academic Press, New York, 31\textendash 65. https://doi.org/10.1016/ B978-0-12-597050-1.50006-3.

\bibitem{sch1} Gerben M.~Scheepmaker, Rob M.P.~Goverde, Leo G.~Kroon, 2017, Review of energy-efficient train control and timetabling, {\em European Journal of Operational Research}, {\bf 2017}, 2, 355\textendash376.  DOI: https://doi.org/10.1016/ j.ejor.2016.09.044.

\bibitem{wanp2} Pengling Wang, Rob M.P.~Goverde, 2017,  Multi-train trajectory optimization for energy efficiency and delay recovery on single-track railway lines, {\em Transportation Research Part B}, {\bf 105},  340\textendash 361.  DOI: https://doi.org/10.1016/j.trb.2017.09.012.

\bibitem{wany1} Wang, Y., De Schutter, B., Van den Boom, T. . J., \& Ning, B. (2013). Optimal trajectory planning for trains---A pseudospectral method and a mixed integer programming approach. {\em Transportation Research Part C}, {\bf 29}, 97\textendash 114.

\bibitem{yan1} Songpo Yang, Jianjun Wu, Xin Yang, Huijun Sun, Ziyou Gao, 2018,  Energy-efficient timetable and speed profile optimization with multi-phase speed limits: Theoretical analysis and application. {\em Applied Mathematical Modelling}, {\bf 56}, 32\textendash 50.  DOI:  https://doi.org/10.1016/j.apm.2017.11.017.

\bibitem{yin1} Jiateng Yin, Tao Tang, Lixing Yang, Jing Xun, Yeran Huang, Ziyou Gao, 2017.  Research and development of automatic train operation for railway transportation systems: A survey, {\em Transportation Research Part C}, {\bf 85}, 548\textendash 572. DOI: 10.1016/j.trc.2017.09.009.

\bibitem{zhu3} X Zhuan and X Xia, 2008, Speed regulation with measured output feedback in the control of heavy haul trains, \textit{Automatica}, {\bf 44}, 242--247.

\end{thebibliography}
\end{document}